\documentclass[final,3p,times]{elsarticle}

\usepackage{amssymb}
\usepackage{graphicx,setspace,color}
\usepackage{psfrag}
\usepackage{amsfonts}
\expandafter\let\csname equation*\endcsname\relax
\expandafter\let\csname endequation*\endcsname\relax
\usepackage{amsmath}
\usepackage{newunicodechar}
\usepackage{amssymb, amsthm}
\usepackage{mathrsfs}
\usepackage{epstopdf}
\usepackage{float}
\usepackage{lineno,hyperref}
\usepackage[capitalise]{cleveref}
\usepackage{color, soul}
\usepackage[normalem]{ulem}
\usepackage{subfig}

\usepackage{dsfont} 
\usepackage{extarrows}
\usepackage{circuitikz}
\usepackage{physics}
\usepackage{bm}

\usepackage{graphicx}
\usepackage{caption}
\usepackage{booktabs}
\usepackage{placeins}

\newunicodechar{●}{$\bullet$}
\newunicodechar{■}{$\blacksquare$}
\newunicodechar{○}{$\circ$}
\newunicodechar{□}{$\square$}
\newcommand{\smallgtrsim}{\mathrel{\scalebox{0.75}{$\boldsymbol{\gtrsim}$}}}

\graphicspath{
  {figures/}
}

\journal{Arxiv}

\begin{document}

\newtheorem{definition}{Definition}[section]
\newtheorem{lemma}{Lemma}[section]
\newtheorem{remark}{Remark}[section]
\newtheorem{theorem}{Theorem}[section]
\newtheorem{proposition}{Proposition}
\newtheorem{assumption}{Assumption}
\newtheorem{example}{Example}
\newtheorem{corollary}{Corollary}[section]
\def\ep{\varepsilon}
\def\Rn{\mathbb{R}^{n}}
\def\Rm{\mathbb{R}^{m}}
\def\E{\mathbb{E}}
\def\hte{\hat\theta}
\def\Schrodinger{$\mathrm{Schr}\ddot{\mathrm{o}}\mathrm{dinger}~$}
\renewcommand{\theequation}{\thesection.\arabic{equation}}

\captionsetup[figure]{justification=raggedright, singlelinecheck=true, format=hang}

\begin{frontmatter}



\title{Bifurcation Structure and Coexisting Dynamical States in Two Inductively Coupled Josephson Junctions}


\author{Ruijie He\fnref{addr1,addr2}}\ead{gravitas_sysu@qq.com}
\author{Shenglan Yuan\corref{cor1}\fnref{addr1}}
\ead{shenglanyuan@gbu.edu.cn}\cortext[cor1]{Corresponding author}

\address[addr1]{\rm Department of Mathematics, School of Sciences, Great Bay University, Dongguan 523000, China }
\address[addr2]{\rm School of Mathematics, Sun Yat-sen University, Guangzhou 510275, China}

\begin{abstract}
Inductively coupled Josephson junctions support coexisting stationary and running states whose selection depends on both circuit parameters and initial conditions.
We investigate the model of a direct current-driven resistively and capacitively shunted junction pair as coupled nonlinear oscillators at zero external magnetic flux, allowing unequal bias, damping and inductive coupling parameters. 
Equilibrium equations, nullclines, and Jacobian spectra are combined with multiple-initial-condition simulations, finite-time Lyapunov spectra, and Poincaré returns to characterize changes in the system dynamics under varying initial conditions and parameter settings.
The fixed-point bifurcation diagrams reveal the creation, destruction or annihilation of equilibria and non-monotonic changes in their number under asymmetric coupling. 
In representative asymmetric parameter slices, nondegenerate saddle-node bifurcations generate or terminate stable equilibria within intervals that also support periodic or chaotic running motion. 
The running dynamics exhibit period-doubling, period-cascading sequences and embedded periodic windows within chaotic regions, while the initial-condition scans resolve their coexistence with stationary states. 
These results connect the local equilibrium structure with the parameter-dependent organization and selection of dynamical regimes in the coupled-junction system.
\end{abstract} 


\end{frontmatter}


\section{Introduction}

A Josephson junction is a superconductive electronic device that consists of two layers of superconductor separated by a thin insulating barrier, which carries a supercurrent at zero voltage.
A direct current superconducting quantum interference device (DC-SQUID) consists of two Josephson junctions connected by a superconducting loop, enabling the system to be sensitive to changes in external magnetic fields.
Its applications include scanning SQUID microscope~\cite{Vu1994ScanningSQUID}, high-resolution scanning~\cite{T_Ishikawa_2006High-Resolution_Magnetic_Image} and biomagnetic measurements~\cite{Mackert1999Magnetoencephalography}.
DC-SQUID is also used in low-noise radio-frequency
amplifiers~\cite{Muck10.1063/1.121490amplifier}, the readout of superconducting flux qubits ~\cite{TANAKA2002300qubit} and neuromorphic circuits~\cite{Mizugaki10.1063/1.108571fluxon}.
The equivalent circuit of a Josephson junction is given by the resistively and capacitively shunted junction (RCSJ) model that consists of an ideal Josephson element $J$ in parallel with a capacitor $C$ and a resistor $R$.
In this work, we assume the two junctions are coupled by the loop inductor $L$, and study how zero-voltage states, running states and chaotic states are organized under varying parameter settings and initial values.

Previous studies have examined deterministic and stochastic dynamics of the two-junction system with different types of coupling. 
Valkering et al. \cite{VALKERING2000137} considered capacitively coupled junctions in the overdamped limit, uncovering non-trivial dynamics including phase locking, devil's staircase and rupture of a toroidal structure.
For resistively coupled models, Kurt et al. \cite{KURT20092229} examined DC-driven models and analyzed the system using bifurcation diagrams, visualization of two-dimensional attractors and Lyapunov spectra. 
They found complex fluctuation patterns of coupling currents, chaotic attractors and a hyperchaotic region where the system admits two positive Lyapunov exponents.
Ben-Jacob and Imry \cite{Ben-Jacob-Dynamics-of-the-DC-SQUID-10.1063/1.328634} studied the current--voltage characteristics and characteristic time scales of an inductively coupled DC SQUID.
Particularly, Baxevanis and Hizanidis~\cite{PhysRevE.111.044214} studied synchronization, multistability and bursting dynamics.
They also established the equilibrium organization and the coexistence of resting and running states for the equal-damping, equal-coupling case with unequal bias, but sampled the parameter plane sparsely: their fixed-point slices use a few values of the bias and their count map is restricted to a narrow window of the asymmetry parameters. 
They also analyzed the fixed-point bifurcations and nullclines at a few representative values of bias and inductance asymmetry.
These studies demonstrate a range of collective behaviors in coupled junctions, while the coupling mechanism and circuit assumptions determine the governing dynamics.  
We therefore examine how parameter asymmetry reorganizes the equilibrium and running dynamics of inductively coupled junctions, resolve the parameter dependence with dense scans, allow unequal damping and unequal coupling, and identify the non-monotonic equilibrium counts and branch reorganizations that these asymmetries produce.

In this work, we consider inductively coupled Josephson junctions with zero external magnetic flux.
In particular, the key normalized system is as follows:
\begin{equation}\label{eq:sde}
\begin{aligned}
\mathrm{d}\phi_1 &= v_1\,\mathrm{d}t, \\
\mathrm{d}v_1 &= \big[i_1 - \beta_1 v_1 - \sin\phi_1 + \kappa_1(\phi_2-\phi_1)\big]\,\mathrm{d}t, \\
\mathrm{d}\phi_2 &= v_2\,\mathrm{d}t, \\
\mathrm{d}v_2 &= \big[i_2 - \beta_2 v_2 - \sin\phi_2 + \kappa_2(\phi_1-\phi_2)\big]\,\mathrm{d}t,
\end{aligned}
\end{equation} 
\nobreak where $i_1, i_2$ are normalized bias currents, $\beta_1, \beta_2$ are damping strengths, and $\kappa_1, \kappa_2$ are coupling strengths.
We focus on three important related questions: (1) how asymmetric junction parameters modify the creation and destruction of equilibria, (2) how the resulting stable branches coexist with running motion, and (3) how initial conditions affect dynamical regimes. 
Addressing these questions requires analysis of equilibria and trajectory-based diagnostics to be interpreted together.
Specifically, we combine calculations of equilibria and Jacobian spectra when one of the parameter varies, multiple-initial-condition simulations, finite-time Lyapunov spectra, Poincaré-section projections and Poincar\'e return periods. 
The equilibrium diagrams reveal non-monotonic changes in the stability and number of fixed points along asymmetric coupling slices. 
Comparisons with the trajectory scans identify intervals in which stable stationary states coexist with running motion, while return period diagnostics resolve period-doubling sequences and periodic windows within the running regimes.

The remainder of this paper is organized as follows.
\cref{sec:coupled-jj-intro} introduces the model of inductively coupled Josephson junctions and derives the governing equations.
 \cref{sec:stability-fixed-points} examines equilibrium stability and the organization of equilibrium branches and nullclines. 
\cref{sec:lyapunov-trajs} presents the Lyapunov-spectrum classifications and representative trajectories. 
\cref{sec:bifurcation_diagram} combines multiple-initial-condition scans with return diagnostics to analyze the stationary and running regimes. 
\cref{sec:conclusion} summarizes the findings and discusses further directions.

\section{Coupled Josephson junctions}\label{sec:coupled-jj-intro}

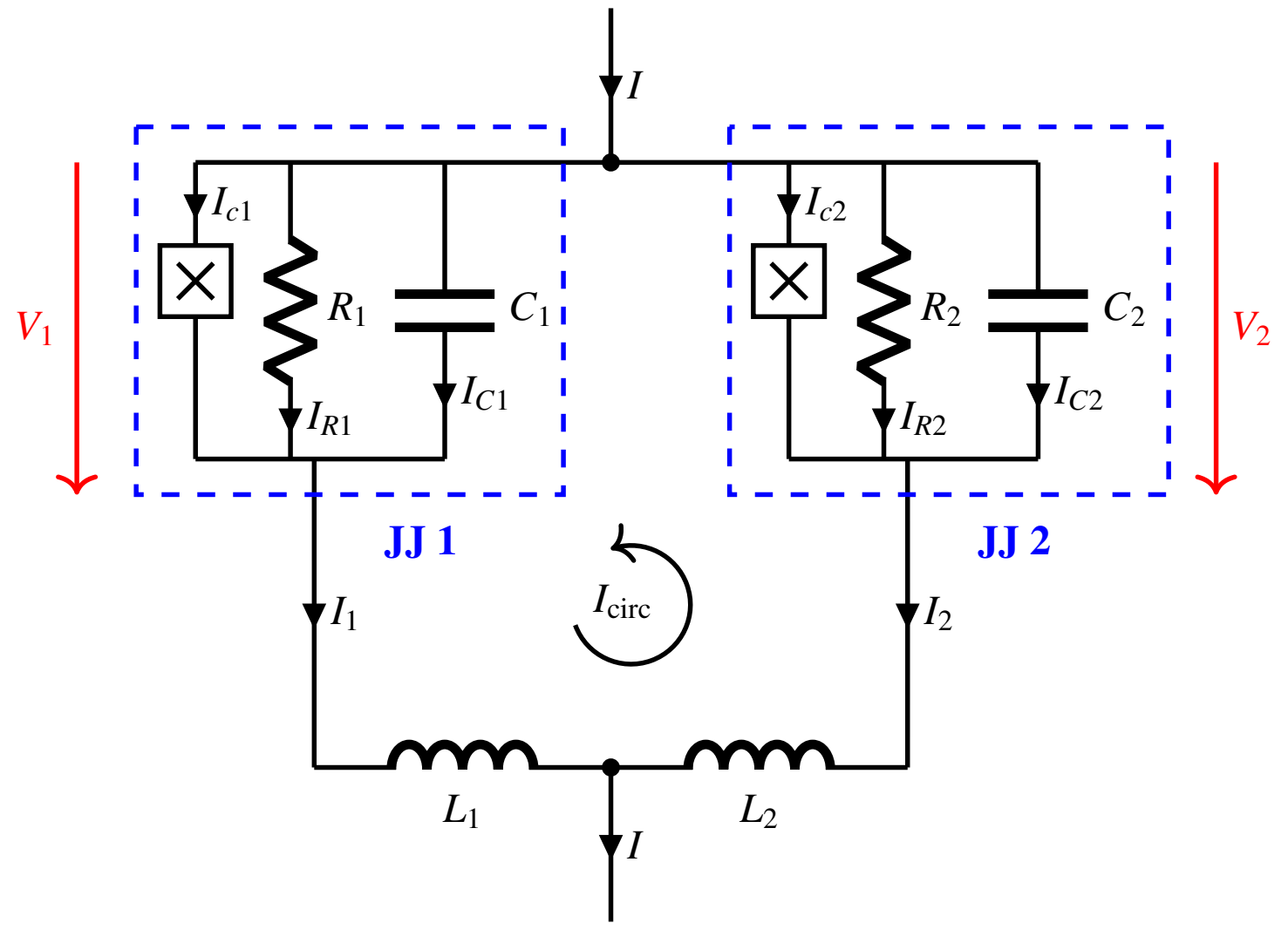
\begin{figure}[htbp]
\centering
\resizebox{0.6\textwidth}{!}{%
\begin{circuitikz}[american, scale=0.7, transform shape]

    \draw[thick] (1.5, 4) -- (8.6, 4);
    \draw[thick] (5, 5.3) to[short, i=$I$, -*] (5, 4);

    \draw[thick] (2.5, -1.1) to[L, l_=$L_1$] (5, -1.1);
    \draw[thick] (5, -1.1) to[L, l_=$L_2$] (7.5, -1.1);
    \draw[thick] (5, -1.1) to[short, i=$I$, *-] (5, -2.4);

    \draw[->, thick] (4.7, 0.1) arc(-160:110:5mm);
    \node at (5.1, 0.3) {$I_{\mathrm{circ}}$};

    \draw[thick] (1.5, 4) to[short, i=$I_{c 1}$] (1.5, 3.3);
    \draw[thick] (1.5, 3.3) -- (1.2, 3.3) -- (1.2, 2.7)
        -- (1.8, 2.7) -- (1.8, 3.3) -- (1.5, 3.3);
    \draw[thick] (1.35, 3.15) -- (1.65, 2.85);
    \draw[thick] (1.35, 2.85) -- (1.65, 3.15);
    \draw[thick] (1.5, 2.7) to[short] (1.5, 1.5);
    \draw[thick] (2.3, 4) to[R, l=$R_1$, i=$I_{R1}$] (2.3, 1.5);
    \draw[thick] (3.6, 4) to[C, l=$C_1$, i=$I_{C1}$] (3.6, 1.5);
    \draw[thick] (1.5, 1.5) -- (3.6, 1.5);
    \draw[thick] (2.5, 1.5) to[short, i=$I_1$] (2.5, -1.1);
    \draw[dashed, thick, blue] (1, 4.3) rectangle (4.6, 1.2);
    \node[blue] at (3.4, 0.8) {\textbf{JJ 1}};
    
    \draw[->, thick, red] (0.5, 4) -- (0.5, 1.2);
    \node[red] at (0.15, 2.6) {$V_1$};

    \draw[thick] (6.5, 4) to[short, i=$I_{c 2}$] (6.5, 3.3);
    \draw[thick] (6.5, 3.3) -- (6.2, 3.3) -- (6.2, 2.7)
        -- (6.8, 2.7) -- (6.8, 3.3) -- (6.5, 3.3);
    \draw[thick] (6.35, 3.15) -- (6.65, 2.85);
    \draw[thick] (6.35, 2.85) -- (6.65, 3.15);
    \draw[thick] (6.5, 2.7) to[short] (6.5, 1.5);
    \draw[thick] (7.3, 4) to[R, l=$R_2$, i=$I_{R2}$] (7.3, 1.5);
    \draw[thick] (8.6, 4) to[C, l=$C_2$, i=$I_{C2}$] (8.6, 1.5);
    \draw[thick] (6.5, 1.5) -- (8.6, 1.5);
    \draw[thick] (7.5, 1.5) to[short, i=$I_2$] (7.5, -1.1);
    \draw[dashed, thick, blue] (6, 4.3) rectangle (9.7, 1.2);
    \node[blue] at (8.4, 0.8) {\textbf{JJ 2}};
    \draw[->, thick, red] (10.1, 4) -- (10.1, 1.2);
    \node[red] at (10.4, 2.6) {$V_2$};
\end{circuitikz}
}
\caption{Equivalent circuit of two coupled junctions  with separate arm
inductances $L_1$ and $L_2$. Each junction is represented by an ideal
Josephson element shunted by a resistor and a capacitor. The total bias
current $I$ enters the common upper node and returns through the midpoint
of the two inductances, i.e.,
$L=L_1+L_2$.}
\label{fig:two_jj_two_inductors}
\end{figure}

We focus on the equivalent circuit of the inductively coupled SQUID in \cref{fig:two_jj_two_inductors}, driven by a constant DC current $I$.
Each junction is represented by a resistively and capacitively shunted junction (RCSJ) model with critical current $I_{c k}$, consisting of a resistance $R_k$ and a capacitance $C_k$. 
The coupling inductance is modeled in the form of two arm inductances $L_1,L_2>0$, without applied external flux. 
Let $\phi_0=\hbar/(2e)$ denote the reduced flux quantum, so that the flux quantum reads $\Phi_0=2\pi\phi_0=h/(2e)$.

Let $\tau$ denote physical time. Here, $I_k$ stands for the total current through branch $k$, and $I_{c k}$ represents the critical current of junction $k$. 
The Josephson currents take the form $I_{c k}\sin\phi_k$, and the Josephson voltage relation is $V_k=\phi_0\,\mathrm{d}\phi_k/\mathrm{d}\tau$.
Applying Kirchhoff's current law, the total current through the $k$-th RCSJ is
\begin{equation}\label{eq:Ij2}
 I_k=I_{c k}\sin\phi_k
     +\frac{\phi_0}{R_k}\frac{\mathrm{d}\phi_k}{\mathrm{d}\tau}
     +\phi_0 C_k\frac{\mathrm{d}^2\phi_k}{\mathrm{d}\tau^2},
 \qquad k=1,2.
\end{equation}
Both branch currents are defined as positive when flowing from the common upper node toward the lower bias-return node as demonstrated in \cref{fig:two_jj_two_inductors}. 
The junction voltages are upper-to-lower voltage drops as in the circuit, and positive counterclockwise circulation $I_{\mathrm{circ}}$ flows downward through junction 1 in the same direction as $I_1$ and upward through junction 2.
Therefore by the law of inductors, RCSJ and Kirchhoff's voltage law, loop constraint yields
\begin{equation}\label{eq:jj-loop}
 V_{L_1} - V_{L_2} = L_1 \frac{\dd I_1}{\dd \tau} - L_2 \frac{\dd I_2}{\dd \tau} \implies
 \phi_0(\phi_2-\phi_1)=L_1I_1-L_2I_2,
\end{equation}
where $V_{L_1}$ and $V_{L_2}$ are voltage differences at the ends of the two inductors.
Using $I = I_1 + I_2$ and letting $L=L_1+L_2$, elementary algebra gives
\begin{equation}\label{eq:jj-branch-currents}
 \begin{aligned}
 I_1&=\frac{L_2}{L}I+\frac{\phi_0}{L}(\phi_2-\phi_1),\\
 I_2&=\frac{L_1}{L}I+\frac{\phi_0}{L}(\phi_1-\phi_2).
 \end{aligned}
\end{equation}
With $I_{\mathrm{circ}}=(I_1-I_2)/2$, the same constraint is
\begin{equation}\label{eq:CC}
 \phi_0(\phi_2-\phi_1)
 =L I_{\mathrm{circ}}+\frac{L_1-L_2}{2}I.
\end{equation}
This constraint adds a bias term to the loop flux~\cite{TescheClarke1977Report} if $L_1 \ne L_2$.
Equal arm inductances yield zero bias, so the loop constraint reduces to $\phi_0(\phi_2-\phi_1)=L I_{\mathrm{circ}}$.

We utilize a device family with equal bare zero-bias plasma frequencies that give the same normalized time:
\begin{equation}\label{eq:jj-shared-time}
 \frac{I_{c1}}{C_1}=\frac{I_{c2}}{C_2}=\phi_0\Omega^2,
 \qquad \Omega>0.
\end{equation}
Given this device family, we set $t=\Omega\tau$ and define
\begin{equation}\label{eq:jj-parameters}
 \beta_k=\frac{1}{\Omega R_k C_k},
 \qquad
 \kappa_k=\frac{\phi_0}{L I_{c k}},
 \qquad 
 i_1 = \frac{L_2 I}{L I_{c1}},
 \qquad 
 i_2 = \frac{L_1 I}{L I_{c2}},
 \qquad k=1,2.
\end{equation}
Here $i_k$ is the normalized constant bias-current component, $\beta_k$ is the dimensionless damping coefficient, and $\kappa_k$ is the dimensionless inductive-coupling strength. 
In particular, the instantaneous total branch current also includes the inductive term:
\begin{equation*}
 \frac{I_k}{I_{c k}}
 =i_k+\kappa_k(\phi_{l}-\phi_k),\qquad (k, l) \in \big\{(1,2), (2, 1)\big\}.
\end{equation*}
The inertial and damping coefficients in normalized time are therefore
$\phi_0 C_k\Omega^2/I_{c k}=1$ and
$\phi_0\Omega/(R_k I_{c k})=\beta_k$, respectively.
Substituting \cref{eq:jj-branch-currents} into \cref{eq:Ij2} and dividing by $I_{c k}$ yields the first-order system
\begin{equation}\label{eq:jj-ode}
 \begin{aligned}
 \dot\phi_1&=v_1,\\
 \dot v_1&=i_1-\beta_1v_1-\sin\phi_1+\kappa_1(\phi_2-\phi_1),\\
 \dot\phi_2&=v_2,\\
 \dot v_2&=i_2-\beta_2v_2-\sin\phi_2+\kappa_2(\phi_1-\phi_2).
 \end{aligned}
\end{equation}
Here $\dot{\square}$ denotes the derivative of variable $\square$ with respect to normalized time $t$, and $v_k=\mathrm{d}\phi_k/\mathrm{d}t = \dot{\phi}_k$.
For fixed components, the single bias links the two drives.
Because the four pairs of parameters satisfy one pair of constraints, the three pairs of normalized parameters $i_k$, $\beta_k$ and $\kappa_k$ are free positive real numbers.

Since this system is invariant with respect to joint translation of $2\pi$, i.e., $(\phi_1, \phi_2) \mapsto (\phi_1+2\pi, \phi_2+2\pi)$ we set $\theta = (\phi_1+\phi_2)/2$ and $\psi = \phi_2 - \phi_1$. 
The equations in \cref{eq:jj-ode} therefore turn into
\begin{equation}\label{eq:coupled_jj-system-transformed}
  \begin{aligned}
     \dot\theta&=\frac{v_1+v_2}{2},\\
     \dot v_1&=i_1-\beta_1v_1-\sin(\theta - \psi/2)+\kappa_1\psi,\\
     \dot\psi&=v_2 - v_1,\\
     \dot v_2&=i_2-\beta_2v_2-\sin(\theta + \psi/2)-\kappa_2\psi,
  \end{aligned}
\end{equation}
in which $\theta \in [0, 2\pi], \psi \in \mathbb{R}$.
On the unwrapped phase space used in system \ref{eq:jj-ode}, this change of variables in system \ref{eq:coupled_jj-system-transformed} is an invertible linear transformation, with inverse mappings $\phi_1=\theta-\psi/2$ and $\phi_2=\theta+\psi/2$.
It gives an equivalent representation of the original dynamics. 
In particular, the Jacobians at corresponding equilibria are similar and have the same spectra and local stability.
The common $2\pi$-translation symmetry further allows $\theta$ to be identified modulo $2\pi$, while $\psi$ remains unwrapped.
Unless explicitly stated otherwise, ``the system'' below refers to \eqref{eq:coupled_jj-system-transformed}, expressed in $(\theta,\psi,v_1,v_2)$ coordinates.
We use "symmetric" to refer to the settings in which two junctions share the same parameters, i.e., $i = i_1 = i_2$, $\beta = \beta_1 = \beta_2$ and $\kappa = \kappa_1 = \kappa_2$.
We also call a parameter setting "asymmetric" if it is not symmetric.

\section{Stability Analysis of Fixed Points}\label{sec:stability-fixed-points}

We calculate the fixed points of system \ref{eq:coupled_jj-system-transformed}.
Setting $\dot{\theta}=\dot{v}_1=\dot{\psi}=\dot{v}_2=0$ gives $v_1=v_2=0$. 
The $v_1$ and $v_2$-nullclines are as follows:
\begin{equation}\label{eq:nullclines_2}
  \begin{aligned}
    i_1-\sin(\theta - \psi/2)+\kappa_1\psi = 0, \\
    i_2-\sin(\theta + \psi/2)-\kappa_2\psi = 0,
  \end{aligned}
\end{equation}
in which $\theta \in [0, 2\pi], \psi \in \mathbb{R}$.
The system Jacobian around the fixed points is 
\begin{equation}\label{Jacobian}
	J = \begin{bmatrix}
		0 &  0 &1/2 & 1/2\\ 
		0 &  0 &-1 & 1\\ 
		-\cos(\theta^* - \frac{\psi^*}{2}) & +\frac12\cos(\theta^*-\frac{\psi^*}{2})+\kappa_1 & -\beta_1 & 0\\
		-\cos(\theta^* + \frac{\psi^*}{2}) & -\frac12\cos(\theta^* + \frac{\psi^*}{2}) - \kappa_2 & 0 & -\beta_2
	\end{bmatrix}.
\end{equation}

There are two settings: either the symmetric setting ($i_1=i_2=i, \kappa_1 = \kappa_2=\kappa, \beta_1 = \beta_2 = \beta$), where all pairs of parameters of the two junctions are the same, or the asymmetric setting.

Consider the symmetric case with $i_1=i_2=i$ and $\kappa_1 = \kappa_2 = \kappa$.
Since $v_1=v_2=0$ at equilibrium, the equilibrium equations do not depend on $\beta_1$ or $\beta_2$.
Adding and subtracting these equations gives
\begin{equation}\label{eq:nullclines_sym}
  \begin{aligned}
    \sin\theta \cos(\psi/2) &= i, \\
    \cos\theta\sin(\psi/2) &= -\kappa\psi.
  \end{aligned}
\end{equation}
\nobreak{}When the two junctions are fully synchronized, i.e., $\psi \equiv 0$, the system \eqref{eq:jj-ode} falls back to the single-junction case.
Note that $\sin\theta \cos(\psi/2) = i$ simplifies to $\sin \theta = i$, as $i$ denotes the non-negative normalized current. 
The fully-synchronized system admits fixed points only when $i \leq 1$.
At $i=1$ this system has exactly one fixed point at $(\theta, \psi) = (\pi/2, 0)$, and the characteristic polynomial of the Jacobian $J$ in \cref{Jacobian} is
\begin{equation}\label{eq:char_poly}
  \lambda (\beta + \lambda) (\lambda^2 + \beta\lambda + 2\kappa) = 0, 
\end{equation}
which yields four eigenvalues: 
\[
  \begin{matrix}
    \lambda_1 = 0, & \lambda_2 = -\beta 
    & \displaystyle \lambda_3 = -\frac{\beta - \sqrt{\beta^2 - 8 \kappa}}{2}, 
    & \displaystyle \lambda_4 = -\frac{\beta + \sqrt{\beta^2 - 8 \kappa}}{2}.
  \end{matrix}
\]
For $\beta^2 \geq 8\kappa$, the $\lambda_3$ and $\lambda_4$ are negative real numbers; 
when $\beta^2 < 8\kappa$, those two eigenvalues form a complex-conjugate pair with negative real parts $-\beta/2$.
This shows that $(\theta, \psi, v_1, v_2) = (\pi/2, 0, 0, 0)$ is a non-hyperbolic equilibrium whose linearization is marginally stable.

Using the parameters set above (symmetric setting with $i_1 = i_2 = 1$) and letting $\vartheta = \theta - \pi/2$, $s = (v_1 + v_2) / 2$ and $d = v_2 - v_1$ yield the following transformed system on this very occasion:
\begin{equation}\label{eq:nullclines_sym_usd_transformed}
  \begin{aligned}
    \dot{\vartheta} &= s, \\
    \dot{s} &= 1 -\beta s - \cos \vartheta \cos \frac{\psi}{2}, \\
    \dot{\psi} &= d, \\
    \dot{d} &= -\beta d + 2\sin \vartheta \sin \frac{\psi}{2} - 2\kappa \psi.
  \end{aligned}
\end{equation}
The fixed point $(\theta, \psi, v_1, v_2) = (\pi/2, 0, 0, 0)$ corresponds to $(\vartheta, \psi, s, d) = (0, 0, 0, 0)$ in the new coordinates used in system \ref{eq:nullclines_sym_usd_transformed}.
The Jacobian of this system at $(\vartheta, \psi, s, d) = (0, 0, 0, 0)$ shares the same set of eigenvalues as that of system \ref{eq:coupled_jj-system-transformed} at $(\theta, \psi, v_1, v_2) = (\pi/2, 0, 0, 0)$.
Moreover, the two-dimensional synchronization manifold $\mathcal{S} = \{(\vartheta, \psi, s, d): \psi = d = 0\}$ is invariant, on which the system dynamics is reduced to
\begin{equation}\label{eq:reduced_dyn_synchronization_manifold_i1}
  \begin{aligned}
    \dot{\vartheta} &= s, \\
    \dot{s} &= 1 -\beta s - \cos \vartheta,
  \end{aligned}
\end{equation}
with eigenvalues $\lambda_1 = 0, \lambda_2 = -\beta$, indicating a one-dimensional smooth center manifold \cite{VANDERBAUWHEDE1987209}. 

Therefore, we let $s = h(\vartheta)$ be the local center manifold, where $h \in C^\infty\big((-\delta, \delta), \mathbb{R}\big)$ satisfying $h(0) = 0$ and $h'(0) = 0$ for some $\delta > 0$.
On the one hand, since the center manifold is invariant, we have $\dot{s} = h'(\vartheta)\dot{\vartheta} = h'(\vartheta)h(\vartheta)$.
On the other hand, according to system \ref{eq:reduced_dyn_synchronization_manifold_i1}, one gets $\dot{s} = 1- \beta s - \cos \vartheta$.
These constraints give the equation $h'(\vartheta)h(\vartheta) = 1- \beta h(\vartheta) - \cos \vartheta$.
Expanding $\cos \vartheta$ to fourth order and writing $h(\vartheta)$ to be $h(\vartheta) = a \vartheta^2 + b \vartheta^3 + O(\vartheta^4)$, we have
\begin{equation}
  \begin{aligned}
    2a^2 \vartheta^3 &= - \beta \left(a \vartheta^2 - \frac{\vartheta^2}{2\beta} + b \vartheta^3\right) + O(\vartheta^4),
  \end{aligned}
\end{equation}
which gives $a = 1/(2\beta)$ and $b = -1/(2\beta^3)$.
Thus, the center manifold satisfies $s = h(\vartheta) = \vartheta^2/(2\beta) - \vartheta^3/(2\beta^3) + O(\vartheta^4) \allowbreak= \vartheta^2/(2\beta) + O(\vartheta^3)$.
There exists a sufficiently small $\rho$, such that $h(\vartheta) \ge \vartheta^2/(4\beta) > 0$ for all $0 < |\vartheta| < \rho$.
This tells that when $\vartheta < 0$, the trajectory approaches origin from the left; when $\vartheta > 0$ the trajectory travels away from the origin, and eventually reaches $\rho$ in finite time.
Therefore, we conclude that nearby trajectories in some neighborhood of the non-hyperbolic equilibrium $(\theta, \psi, v_1, v_2) = (\pi/2, 0, 0, 0)$ are attracted transversely toward the center manifold $s = h(\vartheta)$. 
On this center manifold, the equilibrium is semistable: it attracts trajectories from the $\vartheta < 0$ side and repels them to the $\vartheta > 0$ side.

\section{Fixed-Point Bifurcations and Nullclines}

\begin{figure}[htbp]
  \centering
  \includegraphics[width=\textwidth]{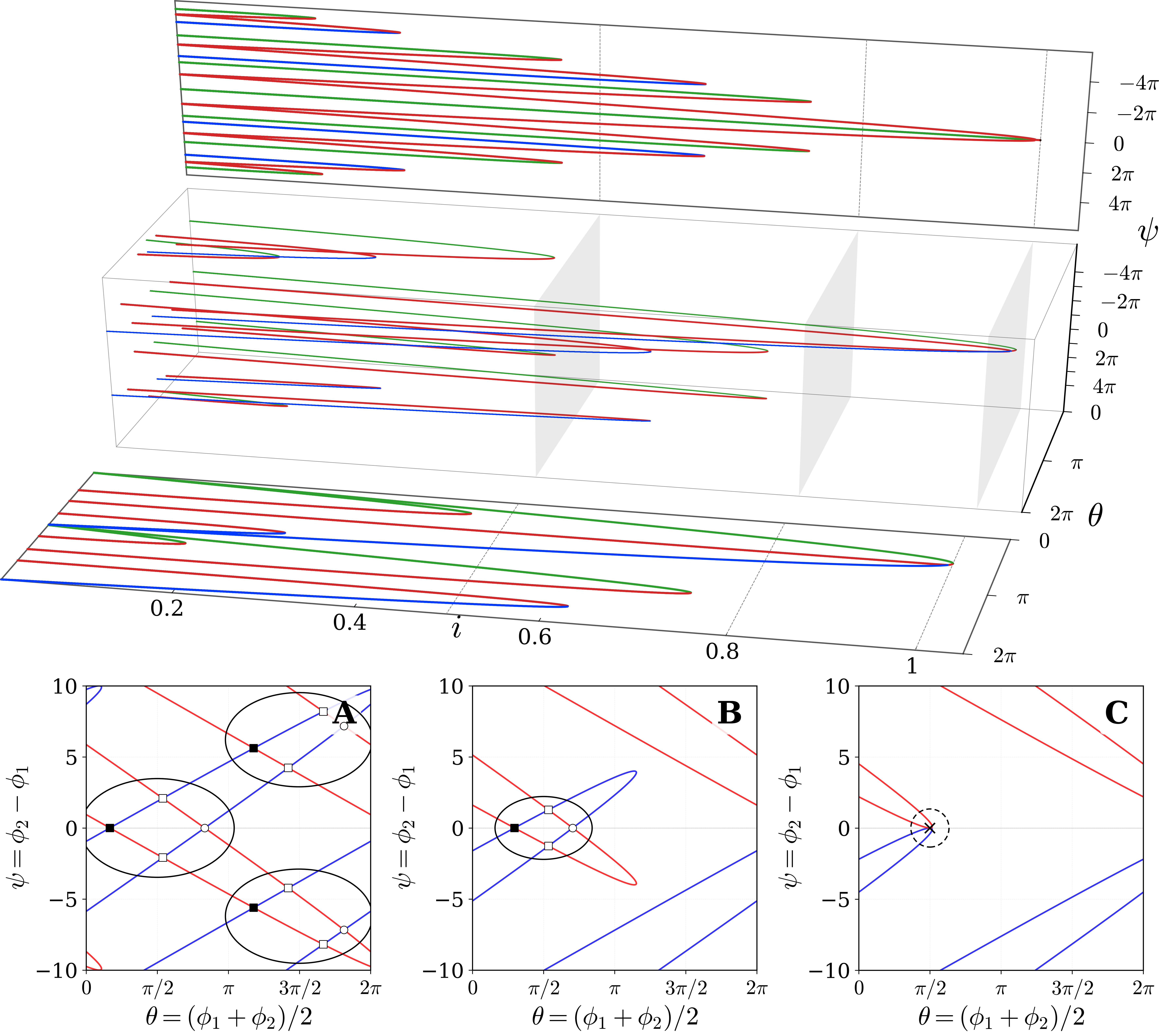}
  \caption{
    Fixed-point bifurcation diagram and nullcline slices of system in the symmetric setting, fixing $\kappa = 0.05$ while varying $i$ from $1.0$ down to $0.01$ with a total of $n = 10{,}000$ frames. 
    The top three panels are the 3D bifurcation diagram over $i$ and its two projections, the $i$--$\psi$ projection on top and the $i$--$\theta$ projection below it.
    Three 2D slices at the bottom labeled (A)--(C) are nullcline slices at $i = 0.5, 0.8, 1.0$, each corresponding to the shaded slice marked in the 3D bifurcation diagram.
    In the 3D bifurcation diagram and its projections, green, red and blue curves trace the stable foci, the saddle foci and the saddle nodes, respectively; in the nullcline slices (A)--(C), red and blue curves correspond to the $v_2$ and $v_1$ nullclines. 
    Filled squares, empty circles and empty squares stand for stable foci, saddle nodes and saddle foci, respectively.
    At $i=1$, a black diamond in the 3D diagram, black dots in its projections and a black cross in (C) mark the non-hyperbolic equilibrium.
  }\label{fig:sym_i_nullclines}
\end{figure}

\begin{figure}[htbp]
  \centering
  \includegraphics[width=\textwidth]{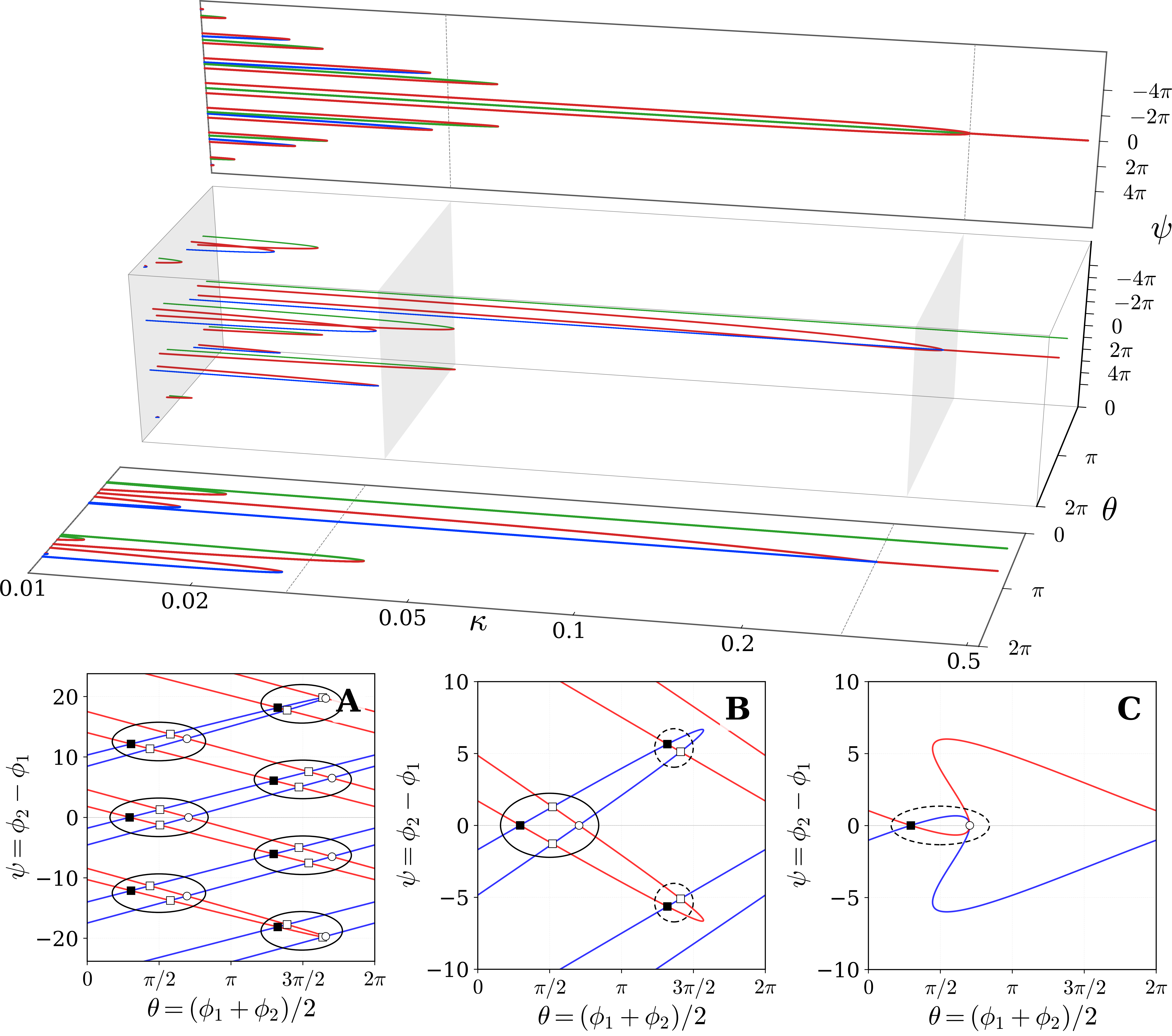}
  \caption{
    Fixed-point bifurcation diagram and nullcline slices in the symmetric setting, fixing $i = 0.8$ while varying $\kappa$ log-spaced from $0.01$ to $0.5$ with a total of $n = 10{,}000$ frames. 
    The nullcline slices (A)--(C) are taken at $\kappa=0.01,0.03,0.3$. 
    The panel arrangement and visual conventions are the same as \cref{fig:sym_i_nullclines}.
  }\label{fig:sym_kappa_nullclines}
\end{figure}

\begin{figure}[htbp]
  \centering
  \includegraphics[width=\textwidth]{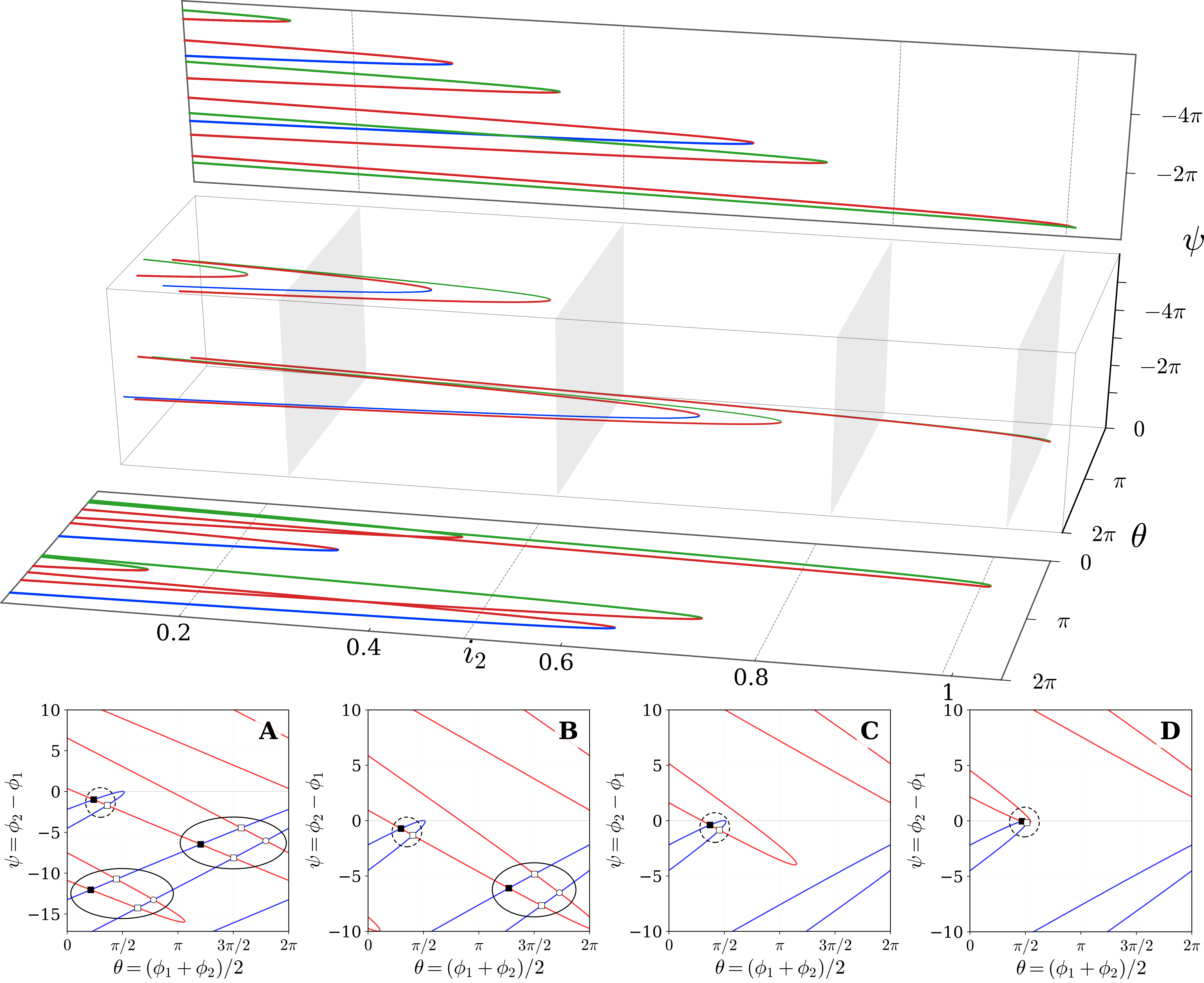}
  \caption{
    Fixed-point bifurcation diagram of the system and nullcline slices in the asymmetric setting, fixing $i_1 = 1$ and $\kappa_1 = \kappa_2 = 0.05$ while varying $i_2$ from $1.0$ down to $0.01$, with a total of $n = 10{,}000$ log-space-distributed frames. 
    The nullcline slices (A)--(D) are taken at $i_2=0.2,0.5,0.8,0.99$. 
    Colors, markers and panel layout inherit those in \cref{fig:sym_i_nullclines}.
  }\label{fig:asym_i2_nullclines}
\end{figure}

\begin{figure}[htbp]
  \centering
  \includegraphics[width=\textwidth]{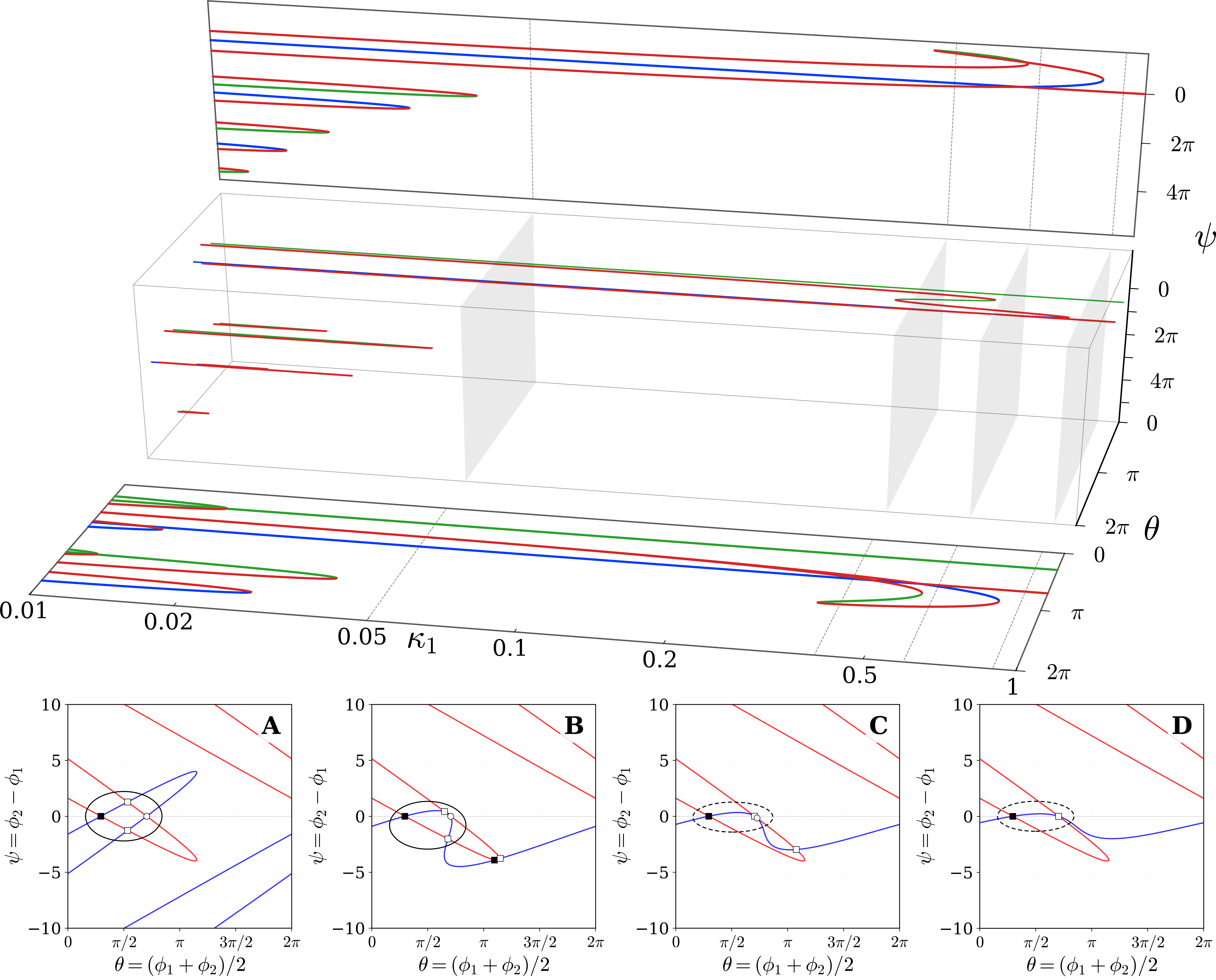}
  \caption{
    Fixed-point bifurcation diagram of the system and nullcline slices in the asymmetric setting, fixing $i_1 = i_2 = 0.8$ and $\kappa_2 = 0.05$ while varying $\kappa_1$ from $0.01$ to $1.0$, with a total of $n = 10{,}000$ log-space-distributed frames. 
    The nullcline slices (A)--(D) are taken at $\kappa_1=0.05,0.4,0.6,0.9$. Colors, markers and panel layout inherit those in \cref{fig:sym_i_nullclines}.
  }\label{fig:asym_kappa1_nullclines}
\end{figure}

To characterize the fixed-point bifurcations, we construct the 3D fixed-point bifurcation diagrams under predetermined values of $i$ and $\kappa$ by keeping one of the parameters fixed and varying the other one. 
Particularly, we consider a parameter interval with other parameters fixed, find the intersections of $v_1$ and $v_2$ nullclines in \cref{eq:nullclines_2} at each parameter step in the interval and stack the results. 
For all figures in this section, we fix $\beta = \beta_1 = \beta_2 = 0.1$.
A hyperbolic equilibrium is called a stable focus when every eigenvalue has a negative real part and at least one pair is complex, a saddle node when its spectrum is real with both positive and negative eigenvalues, and a saddle focus when it has one positive real eigenvalue together with a complex-conjugate pair. 
This yields the 3D fixed-point bifurcation diagram presented in~\cref{fig:sym_i_nullclines,fig:sym_kappa_nullclines,fig:asym_i2_nullclines,fig:asym_kappa1_nullclines}, in which both 2D projections and 10,000 sampled parameter values are provided to give a multi-perspective view of the fixed-point bifurcations.

\cref{fig:sym_i_nullclines} illustrates the fixed-point bifurcations of system in the symmetric setting, fixing $\kappa = 0.05$ while varying $i$ from $1.0$ down to $0.01$.
We sample values of $i$ and plot the resulting fixed points in the top 3D bifurcation diagram, in which blue, red and green curves correspond to the trajectories of saddle nodes, saddle foci and stable foci, respectively. 
We choose $i = 0.5$, $0.8$, $1.0$ from the slices and present the 2D slices (A)--(C) below with classified fixed points and $v_1$, $v_2$ nullclines, in which stable foci, saddle nodes, saddle foci are marked by filled squares, empty circles, empty squares, respectively.
In this symmetric setting, the $v_1$ nullcline is a mirror reflection of the $v_2$ nullcline about $\psi = 0$, and therefore the number of fixed points with $\psi \neq 0$ shall be even. 
Since $i_1$ and $i_2$ are external bias currents, changing $i$ shifts the nullcline curves. 
The 3D bifurcation diagram depicts a saddle-saddle (SS) bifurcation: the two saddle foci with $\psi \ne 0$ in (B) merge with the saddle node at $\psi = 0$, which persists.
Then the remaining equilibria, a saddle focus and a stable focus, collide, forming the saddle-node (SN) bifurcation at $i = 1$, resulting in a non-hyperbolic critical point at $(\theta, \psi) = (\pi/2, 0)$, which disappears as $i$ increases further.
In the $i$--$\psi$ projection, the bifurcations occur in pairs for smaller $i$: two SN bifurcations followed by two SS bifurcations.
The same pattern is visible in the nullcline slices (A) and (B). 
Here a group means the four fixed points created, destroyed or annihilated in pairs as $i$ decreases.
In symmetric settings, the group off the synchronization line $\psi = 0$ always has a mirror-image about this line. 
Conversely, as $i$ increases toward $1$, each group loses two fixed points in an SN bifurcation after losing a saddle node and a saddle focus in a SS bifurcation. 
In particular, when $i$ increases, the saddle focus with the larger $|\psi|$ and the saddle node annihilate, leaving one saddle focus and one stable focus in each group. 
Then the remaining fixed points of the two groups continue to collide and annihilate as $i$ gets larger.
In a bird's-eye view, the $\theta$ coordinate of saddle-node bifurcations advances monotonically when $i$ increases from $0.01$.
With $\theta$ unwrapped across the $2\pi$ boundary, the bifurcation positions exhibit some periodic pattern: each bifurcation pair occurs at the same $i$ and wraps back to a small $\theta$ value when $\theta_*$ crosses $2\pi$. 
The annihilating pairs alternate between (stable focus, saddle focus) pairs and (saddle focus, saddle node) pairs, so that the pattern of the bifurcation positions repeats approximately every three consecutive pairs of bifurcations.

\cref{fig:sym_kappa_nullclines} shares the same setting and display conventions as \cref{fig:sym_i_nullclines}, with $i = 0.8$ fixed and $\kappa$ log-spaced from $0.01$ to $0.5$. 
The nullclines at $\kappa = 0.01, 0.03, 0.3$ are plotted in (A)--(C). 
As $\kappa$ increases, the fixed points annihilate in pairs, and this pattern repeats as $\kappa$ grows further. 
The nullclines become progressively flatter, and at $\kappa \approx 0.299$ the two saddle foci at $\theta \approx 0.7\pi$ merge with the saddle node on $\psi = 0$, which persists after this collision. 
For $\kappa \smallgtrsim 0.3$, a single stable focus at $\theta \approx 0.3\pi$ and one saddle at $\theta \approx 0.7\pi$ remain.

\begin{figure}[htbp]
    \centering
    \includegraphics[width=0.8\textwidth]{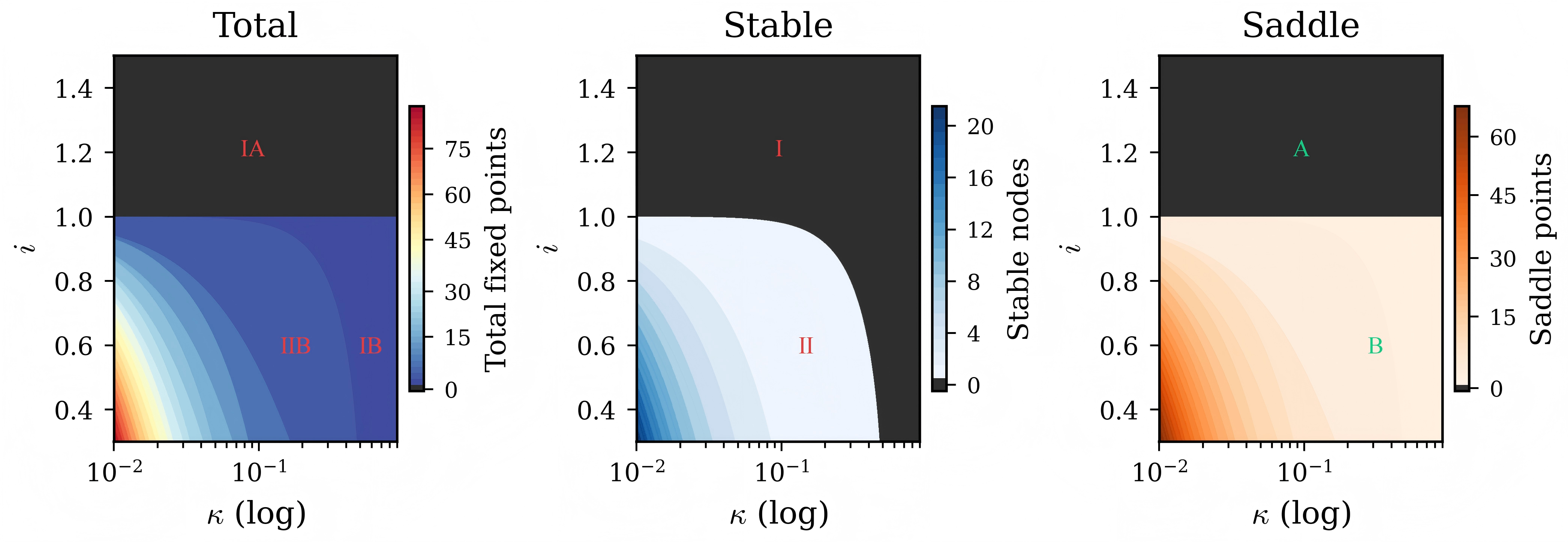}
    \caption{
      Number of fixed points of the system in the symmetric case $i_1=i_2=i$, $\kappa_1=\kappa_2=\kappa$ over $(\kappa,i)\in[0.01,0.9]\times[0.3,1.5]$ with log-scale axis on $\kappa$. 
      The left panel depicts the total number of fixed points, the middle shows that of stable fixed points and the rightmost one reveals the number of saddles. 
      Connected regions in each panel, colored according to the color map on the immediate right side of the panel, are labeled with Roman numerals and capitalized letters.
      Zero counts are encoded in black and higher non-zero numbers of fixed points are filled with color that is closer to red (left) or darker (middle, right).
    }
    \label{fig:num_fix_sym}
\end{figure}

\cref{fig:num_fix_sym} condenses the fixed-point structure of the symmetric case into global count maps over the $(\kappa, i)$ parameter plane.
The horizontal line $i=1$ cleanly separates region~I ($i > 1$, no fixed points) from other regions. 
At $i=1$ the unique equilibrium $(\theta, \psi, v_1, v_2) = (\pi/2, 0, 0, 0)$ has one zero eigenvalue and three eigenvalues that have negative real parts.
As either $i$ or $\kappa$ decreases, $N_{\text{tot}}$ grows monotonically; in the limit $\kappa \to 0$ the nullclines oscillate with undiminished amplitude, producing infinitely many intersections.
The middle and right panels decompose the total count into stable and saddle contributions, respectively, revealing that saddles proliferate roughly three times faster than stable points.
This is consistent with the four-fixed-point groups identified in \cref{fig:sym_i_nullclines}.
The vertical line $\kappa=0.05$ and the horizontal line $i=0.8$ marked on the map correspond to the slice paths in panels (A)--(C) of \cref{fig:sym_i_nullclines} and of \cref{fig:sym_kappa_nullclines}, respectively; the latter ends in region~II with $N_{\text{tot}} = 2$.

\cref{fig:asym_i2_nullclines} plots the fixed-point bifurcation diagram of the system in the asymmetric setting, with $i_1 = 1$, $\kappa_1 = \kappa_2 = 0.05$ fixed and $i_2$ log-spaced from $1.0$ down to $0.01$. 
Panels (A)--(D) display the nullclines at $i_2 = 0.2, 0.5, 0.8, 0.99$, respectively.
In (A)--(D), the red $v_2$ nullcline moves to the bottom-right direction as $i_2$ decreases, intersecting the blue $v_1$ nullcline, and thus the number of fixed points increases. 
In the $i_2$--$\psi$ projection, the bifurcation positions move toward larger $\psi$. 
Their period pattern is similar to that in \cref{fig:sym_i_nullclines}: each saddle-node bifurcation is followed by a SS bifurcation as $i_2$ decreases. 
In contrast to \cref{fig:sym_i_nullclines}, the $\theta$ coordinate of the bifurcation follows a similar periodic wrapping pattern but with different period.

\cref{fig:asym_kappa1_nullclines} illustrates the fixed-point bifurcations of the system in the asymmetric setting, with $i_1 = i_2 = 0.8$, $\kappa_2 = 0.05$ fixed and $\kappa_1$ log-spaced from $0.01$ to $1.0$. 
The nullcline slices at $\kappa_1 = 0.05, 0.4, 0.6, 0.9$ correspond to panels (A)--(D), respectively.
The nullclines in \cref{fig:asym_kappa1_nullclines} are translated when $i_1$ or $i_2$ changes, and their shape changes with $\kappa_1$: they become flatter as $\kappa_1$ increases.
The change in the shape of the nullclines triggers non-monotonic variations in the stability and number of fixed points even though $\kappa_1$ changes monotonically, as displayed in \cref{fig:asym_kappa1_nullclines} with $\kappa_1 \gtrsim 0.1$, where the $v_1$ and $v_2$ nullclines exhibit notable shape differences. 
When $\kappa_1$ increases from $0.01$ to $0.1$, the fixed-point annihilation follows a pattern similar to that in \cref{fig:sym_i_nullclines,fig:asym_i2_nullclines}. 
A saddle-node bifurcation is always followed by an SS bifurcation, then another saddle-node bifurcation, and so on.
However, for $\kappa_1 > 0.1$, a saddle-node bifurcation at $\kappa_1 \approx 0.36$ creates a stable focus and a saddle focus, both with $\theta$ slightly above $\pi$. At $\kappa_1 = 0.4$, the two new fixed points are created at $(\theta, \psi) = (1.1\pi, -1.25\pi)$ and $(1.15\pi, -1.2\pi)$, respectively, as illustrated in (B). 
As $\kappa_1$ increases further, the newborn stable focus migrates toward smaller $\theta$, and annihilates with the saddle focus near $(\theta, \psi) \approx (0.87\pi, -0.94\pi)$ at $\kappa_1 \approx 0.57$, as illustrated in (C), not long after a transcritical bifurcation between two saddle equilibria near $\psi = 0$ at $\kappa_1 \approx 0.55$.
At $\kappa_1 = 0.6$, only four fixed points remain: a stable focus at $(0.3\pi, 0)$, a saddle focus at $(0.7\pi, 0)$, a saddle node at $(0.73\pi, -0.05\pi)$, and a saddle focus at $(1.08\pi, -0.95\pi)$.
Finally, the saddle node at $(\theta, \psi) = (0.89\pi, -0.46\pi)$ and the saddle focus at $(0.9\pi, -0.49\pi)$ annihilate at $\kappa_1 \approx 0.814$. 
The two remaining fixed points are a stable focus at $(0.3\pi, 0)$ and a saddle focus at $(0.7\pi, 0)$, both with $\psi = 0$, as depicted in panel (D).

\begin{figure}[htbp]
    \centering
    \includegraphics[width=0.8\textwidth]{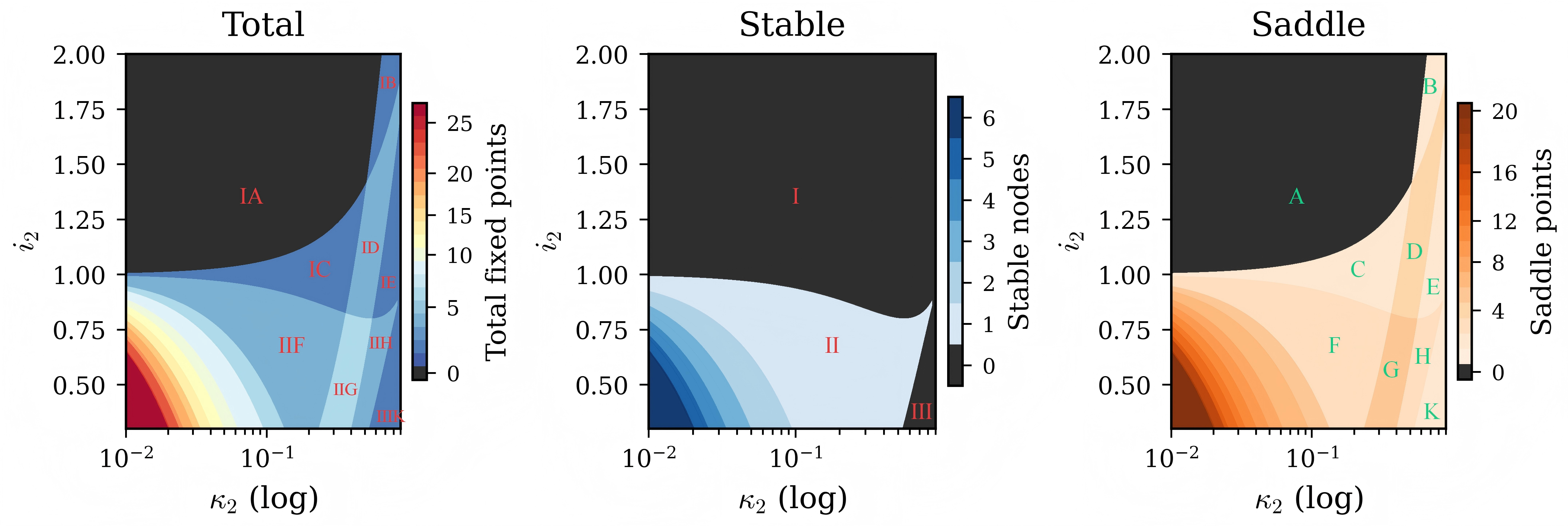}
    \caption{
    Number of fixed points of the system in the asymmetric case under $i_1=0.8$, $\kappa_1=0.05$ over $(i_2,\kappa_2)\in[0.3,2]\times[0.01,0.9]$ with log-scale axis on $\kappa_2$. 
    Panel organization, color scales and region labels remain the same as \cref{fig:num_fix_sym}.
  }
    \label{fig:num_fix_asym}
\end{figure}

\cref{fig:num_fix_asym} sketches the fixed-point count map of the system for the asymmetric case with $i_1 = 0.8$ and $\kappa_1 = 0.05$ fixed, using the same color scheme as \cref{fig:num_fix_sym}.
The critical boundary $i=1$ is deformed into an $i_2$-dependent curve, which rises with $\kappa_2$ and intersects a narrow vertical strip on the right that has no counterpart in the symmetric case.
For $i_2 \ll 1.0$ and $\kappa_2 \ll 0.1$, the monotone rule still holds, i.e., smaller $i_2$ or $\kappa_2$ yields more fixed points.  
When $\kappa_2 > 0.1$, the $\dot v_2 = 0$ nullcline differs markedly in shape from the $\dot v_1 = 0$ nullcline, introducing strong nonlinearity. 
To stress the differences and similarities between \cref{fig:num_fix_asym} and \cref{fig:num_fix_sym}, we label regions with Roman numerals and subregions with capital letters. 
Region~I (top) and Region~III (lower right corner) are two regions with zero stable equilibria.
Region~II (at least one stable equilibrium) covers the remaining plane, and sub-regions A, B, C, D, E, F, G, H, K (see the right panel) are marked in the rightmost panel, indicating connected regions of different values of saddle equilibria.
The strip on the right (containing B, D, G, E, H) splits the original structure.
Consequently, fixing $\kappa_2$ and lowering $i_2$, or fixing $i_2 \ll 1.5$ and raising $\kappa_2$ above $0.1$, produces non-monotonic variations in $N_{\text{tot}}$; for example, fixing $i_2 = 0.8$ and increasing $\kappa_2$ from $0.01$ to $0.9$ traverses sub-regions with order II-F $\to$ II-G $\to$ II-H $\to$ III-K, yielding a decrease--increase--decrease sequence in $N_{\text{tot}}$, qualitatively analogous to the saddle cascade in panels (A)--(D) of \cref{fig:asym_kappa1_nullclines}.

\section{Lyapunov Spectrum and Typical Trajectories}\label{sec:lyapunov-trajs}

\paragraph{Boundedness and synchronized dynamics}
For $\beta_1,\beta_2>0$ and $\kappa_1,\kappa_2\geq0$ with $\kappa_1+\kappa_2>0$, the relative phase $\psi=\phi_2-\phi_1$ and both velocities are bounded.
In the variables $(\psi, s, d)$, the system consists of a stable linear part driven by bounded bias and bounded trigonometric terms.
Since $\dot\psi=v_2-v_1$, the difference between the time-averaged velocities satisfies $T^{-1}\int_0^T(v_2-v_1)\,\mathrm{d}t=[\psi(T)-\psi(0)]/T$ that vanishes as $T \to \infty$, thus, the difference between the mean voltages vanishes as the averaging time increases.

For symmetric couplings $\kappa_1=\kappa_2=\kappa$, the tilted mechanical energy $E=\tfrac12(v_1^2+v_2^2)+\sum_{j=1}^2(1-\cos\phi_j-i_j\phi_j)+\tfrac\kappa2\psi^2$ satisfies $\dot E=-\beta_1v_1^2-\beta_2v_2^2\leq0$. 
For trapped trajectories with bounded unwrapped phases, LaSalle's invariance principle~\cite{khalil2002nonlinear} gives convergence to the set of fixed points, which is the largest invariant subset of $\{v_1=v_2=0\}$. 

In the fully symmetric case, the synchronization manifold $\mathcal S$ introduced in \cref{sec:stability-fixed-points} remains invariant for arbitrary common drive $i$.
On $\mathcal S$, the coupling vanishes, and the system reduces to \cref{eq:reduced_dyn_synchronization_manifold_i1} with constant $1$ replaced by the common drive $i$.
The resulting two-dimensional autonomous flow has no chaotic attractors, but exhibits trapped equilibria or periodic running states.

\begin{figure}[htbp]
  \centering
  \includegraphics[width=0.6\textwidth]{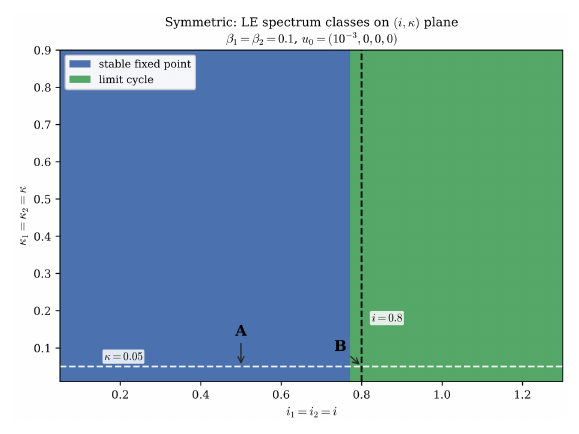}
  \caption{
      Finite-time Lyapunov-spectrum classification on the symmetric $(i,\kappa)$ parameter plane, with $\beta_1=\beta_2=0.1$ and initial state $(\theta, \psi, v_1, v_2) = (5 \times 10^{-4}, -10^{-3}, 0, 0)$. 
      Blue and green denote fixed points and limit cycles, respectively, as observed in the limit dynamics. 
      The marked cases A and B are displayed in \cref{fig:sym_le_traj_ab}.}
  \label{fig:sym_le_spec_trajs}
\end{figure}

\paragraph{Numerical Scanning of Lyapunov Spectrum}
We then further explore different types of dynamics by calculating their Lyapunov exponents using the Benettin--QR \cite{benettin1980lyapunov} method with an adaptive solver \texttt{AutoTsit5(Rodas5P())}.
Unless stated otherwise, we set transient time $T_{\mathrm{tr}}=300$ units and measurement time $T_{\mathrm{me}}=1500$ units, and the initial-condition maps use $T_{\mathrm{tr}} = 100$, $T_{\mathrm{me}} = 400$.
QR reorthonormalization is performed every unit time.
Let $\lambda_1, ..., \lambda_4$ denote the estimated four Lyapunov exponents of the system in descending order.
For symmetric settings, we plot the classification of the Lyapunov spectrum with $\beta_1=\beta_2=0.1$, fixing the initial condition $(\theta, \psi, v_1, v_2) = (5 \times 10^{-4}, -10^{-3}, 0, 0)$, and letting the common drive $i_1=i_2=i$ and the coupling $\kappa_1=\kappa_2=\kappa$ vary in $[0.05, 1.3]$ and $[0.01, 0.9]$ respectively, as presented in \cref{fig:sym_le_spec_trajs}. 
We classify the finite-time estimates of Lyapunov exponents into five classes: (a) hyperchaos ($\lambda_1 \ge \lambda_2 > 0$), (b) chaos ($\lambda_1 > \epsilon \ge \lambda_2, \lambda_3, \lambda_4$, and it is not hyperchaos), (c) torus ($\epsilon > \lambda_1 \ge \lambda_2 > -\epsilon > \lambda_3, \lambda_4$), (d) limit cycle ($\epsilon > \lambda_1 > -\epsilon > \lambda_2, \lambda_3, \lambda_4$), and (e) fixed point ($-\epsilon > \lambda_1, \lambda_2, \lambda_3, \lambda_4$).
Most positive estimates of $\lambda_2$ lie within the common tolerance band $\epsilon=0.005$. 
To retain their sign information and spatial distribution in the finite-time classification maps, we use the criterion $\lambda_2>0$ for the hyperchaotic class, which are classified first, and apply $\epsilon=0.005$ only to the remaining classes.

\begin{figure}[htbp]
  \centering
  \includegraphics[width=0.8\textwidth]{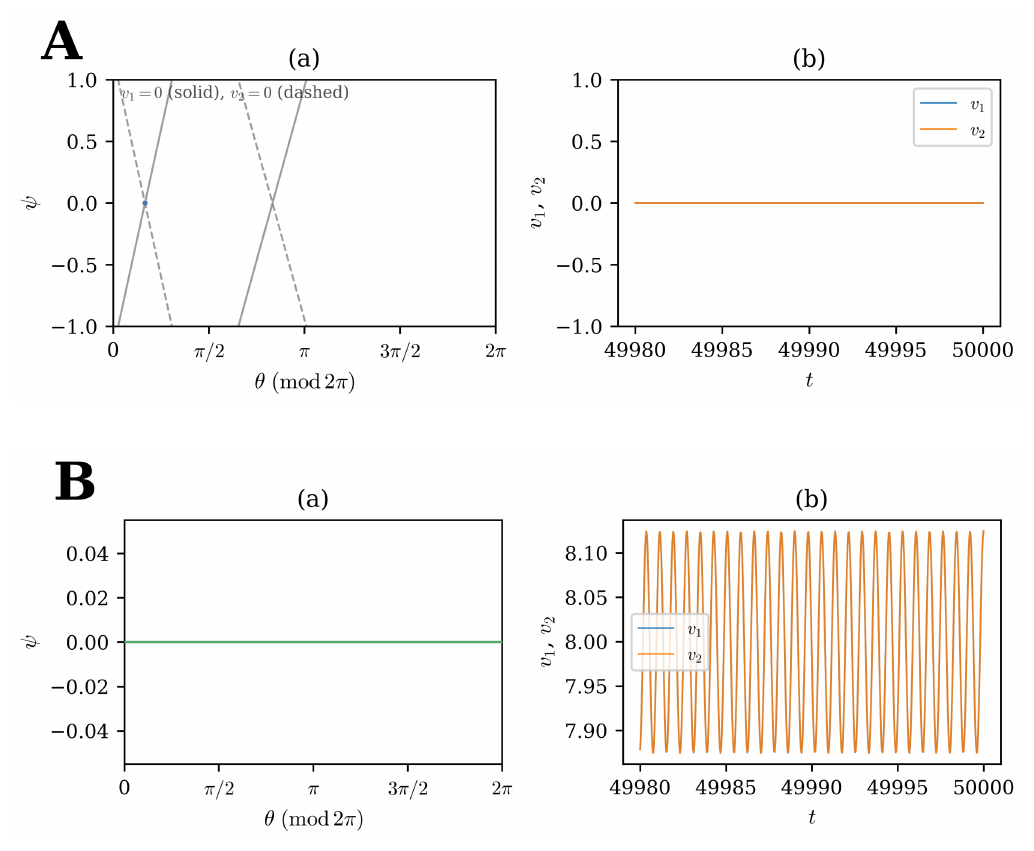}
  \caption{
      Representative trajectories A and B for the symmetric parameter scan in \cref{fig:sym_le_spec_trajs}. 
      Each group displays the last $25{,}000$ time units in the $(\theta\bmod2\pi,\psi)$ projection and the last $20$ time units in the $v_1$ and $v_2$ time series. The initial state of panel (A) is $(\theta, \psi, v_1, v_2) = (0.525, -0.01, 0, 0)$ and that of panel (B) is $(\theta, \psi, v_1, v_2) = (5 \times 10^{-4}, -10^{-3}, 0, 0)$.}
  \label{fig:sym_le_traj_ab}
\end{figure}

For the symmetric map in \cref{fig:sym_le_spec_trajs}, the scanning window $i \times \kappa \in [0.05, 1.3] \times [0.01, 0.9]$ is evenly partitioned into a $52 \times 32$ grid in which blue and green regions correspond to stable fixed points and limit cycles, respectively. 
The fixed-point and periodic-running classes in this scan agree with the regimes of the synchronized reduction.
Along the $\kappa = 0.05$ slice, the boundary between the fixed-point and limit-cycle regions is approximately $i \approx 0.785$. 
We also plot the projected dynamics of the two selected parameter settings (A) and (B) over the last 25,000 time units of a total of 50,000 time units of simulation in \cref{fig:sym_le_traj_ab}. 
The upper setting (A) in \cref{fig:sym_le_traj_ab} belongs to the blue stable-fixed-point region in \cref{fig:sym_le_spec_trajs}. 
The system therefore settles into the stable focus at the left intersection of the $v_1$ and $v_2$ nullclines marked in blue dot.
In the lower setting (B), the system exhibits oscillating $v_1$ and $v_2$ while the phases of the two junctions are locked, resulting in a rotating state. 

\begin{figure}[htbp]
  \centering
  \includegraphics[width=0.7\textwidth]{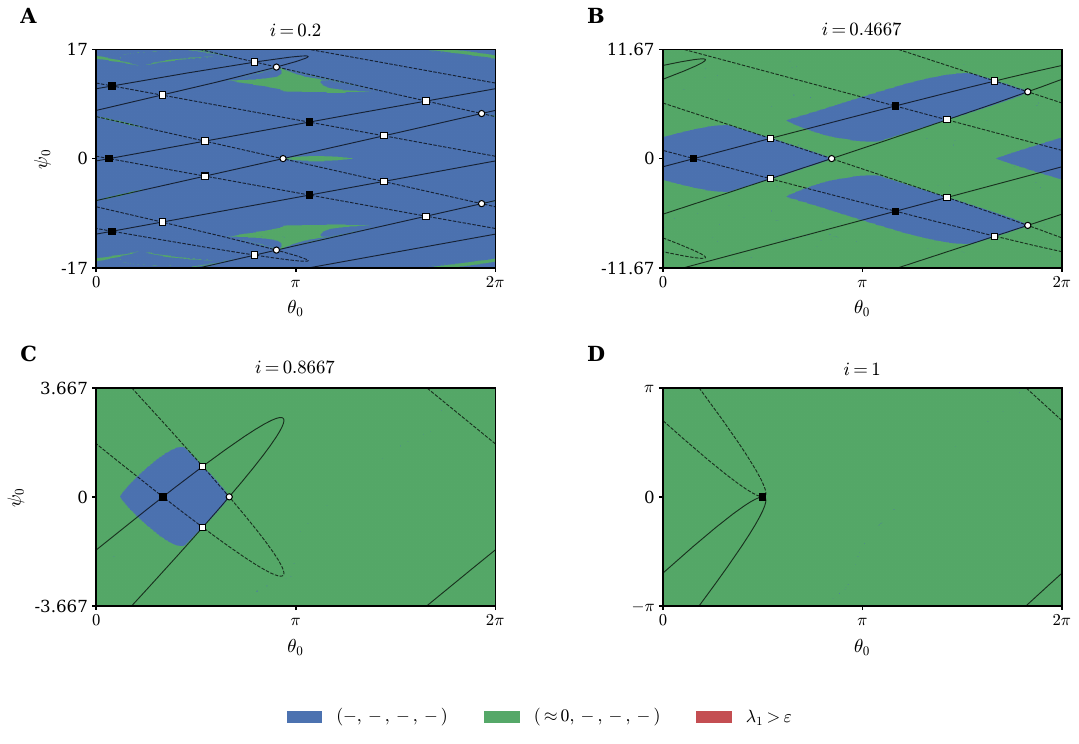}
  \caption{
    Finite-time deterministic Lyapunov classifications while varying symmetric normalized current $i=i_1=i_2$ at fixed $\kappa=\kappa_1=\kappa_2=0.05$. 
    For the symmetric cases $\beta_1=\beta_2=0.1$. 
    The initial-condition plane is $(\theta_0,\psi_0)$ with $v_1(0)=v_2(0)=0$. Each panel uses a $500\times500$ grid. 
    Blue, green and red denote the fixed point, limit cycle and $\lambda_1 > \epsilon$ classes; the isolated red points are transient finite-time positives during synchronization. 
    Solid and dashed curves are the two zero-velocity nullclines, and fixed-point markers follow \cref{fig:sym_i_nullclines}.
  }
  \label{fig:det_le_ic_0}
\end{figure}

\begin{figure}[htbp]
  \centering
  \includegraphics[width=0.8\textwidth]{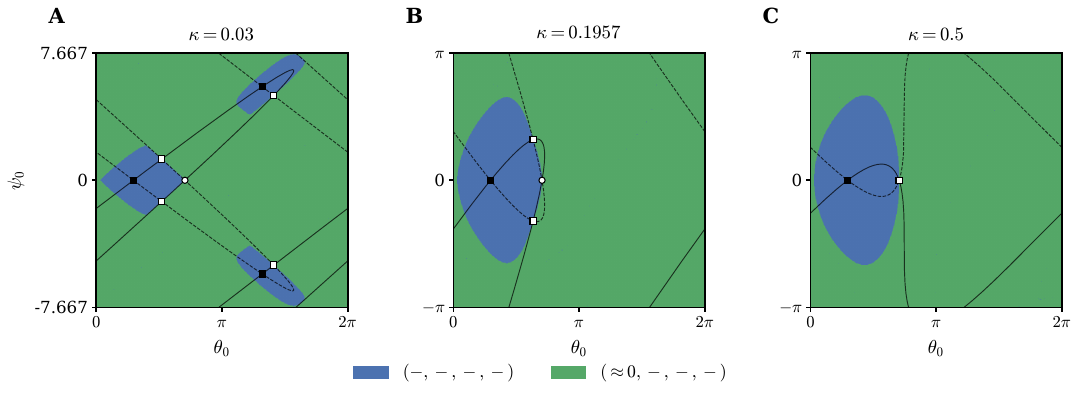}
  \caption{
    Finite-time deterministic Lyapunov classifications while varying symmetric coupling strengths $\kappa$ at fixed $i=0.8$. For the symmetric cases $\beta_1=\beta_2=0.1$. 
    Initial-condition sampling and graphical conventions are as specified in \cref{fig:det_le_ic_0}.
  }
  \label{fig:det_le_ic_1}
\end{figure}

The selected symmetric initial-condition maps in \cref{fig:det_le_ic_0,fig:det_le_ic_1} cover four current values, $i=0.2$, $0.4667$, $0.8667$, and $1.0$, and three coupling values, $\kappa=0.03$, $0.1957$, and $0.5$. 
In \cref{fig:det_le_ic_0}, we fix $\beta_1=\beta_2=0.1$, $\kappa_1 = \kappa_2 =0.05$ and increase $i = i_1 = i_2$. 
The blue regions contract as the common drive $i$ increases, while the green region corresponding to limit cycles (or rotating phases) expands.
This forces the blue regions denoting stable foci (or trapped phases) to shrink and deform to separate islands, whose number decreases as the number of stable foci decreases with increasing $i$.
Moreover, at $i=0.4667$ and $i=0.8667$, roughly half of each trapped island's boundary follows the $v_1$ and $v_2$ nullclines.
These boundary segments run from a saddle node to a saddle focus of the same four-fixed-point group depicted in \cref{fig:sym_i_nullclines}. 
Note that there are few red blocks shown in this figure, which are verified as synchronized running state when measuring window is longer.

In \cref{fig:det_le_ic_1}, we fix $i_1 = i_2 =0.8$, $\beta_1=\beta_2=0.1$, and increase $\kappa = \kappa_1 = \kappa_2$.
The first two frames ($\kappa = 0.03$ and $\kappa=0.1957$) display patterns similar to those in the last $3$ frames ($i = 0.4667, 0.8667$ and $i=1$) of \cref{fig:det_le_ic_0}. 
As $\kappa$ grows, the nullclines become less fluctuating, producing fewer stable foci: the multiple blue regions when $\kappa=0.03$ pave the way to a single broad region in the displayed windows at $\kappa=0.1957$ and $0.5$.

Notably, once the three saddle fixed points in the last group of four fixed points at $\kappa=0.1957$ collide to form a single saddle focus at $\kappa=0.5$, the edge of the blue region no longer aligns with the $v_1$ and $v_2$ nullclines.

\begin{figure}[htbp]
  \centering
  \includegraphics[width=0.6\textwidth]{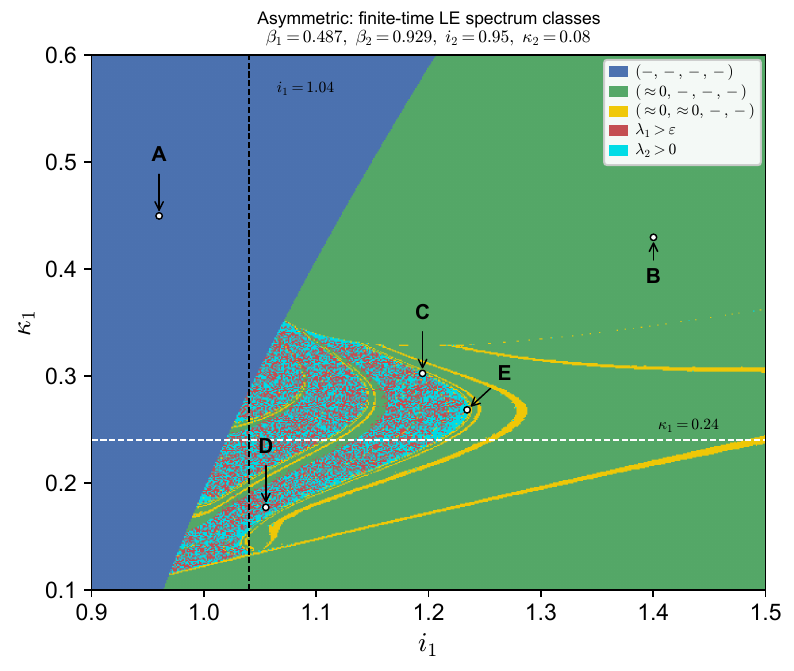}
  \caption{
    Finite-time Lyapunov-spectrum classification on a $500\times500$ grid in the asymmetric $(i_1,\kappa_1)$ plane, with $(\beta_1,\beta_2,\allowbreak i_2,\kappa_2)\allowbreak=(0.487,0.929,\allowbreak0.95,0.08)$ and zero initial state. Blue, green, yellow, red and cyan indicate fixed point, limit cycle, torus, chaos and hyperchaos, respectively, as observed limit dynamics.
    Representative points A to E correspond to \cref{fig:asym_le_traj_ab,fig:asym_le_traj_c,fig:asym_le_traj_d,fig:asym_le_traj_e}.
  }
  \label{fig:asym_le_spec_trajs}
\end{figure}

\begin{figure}[htbp]
  \centering
  \includegraphics[width=0.8\textwidth]{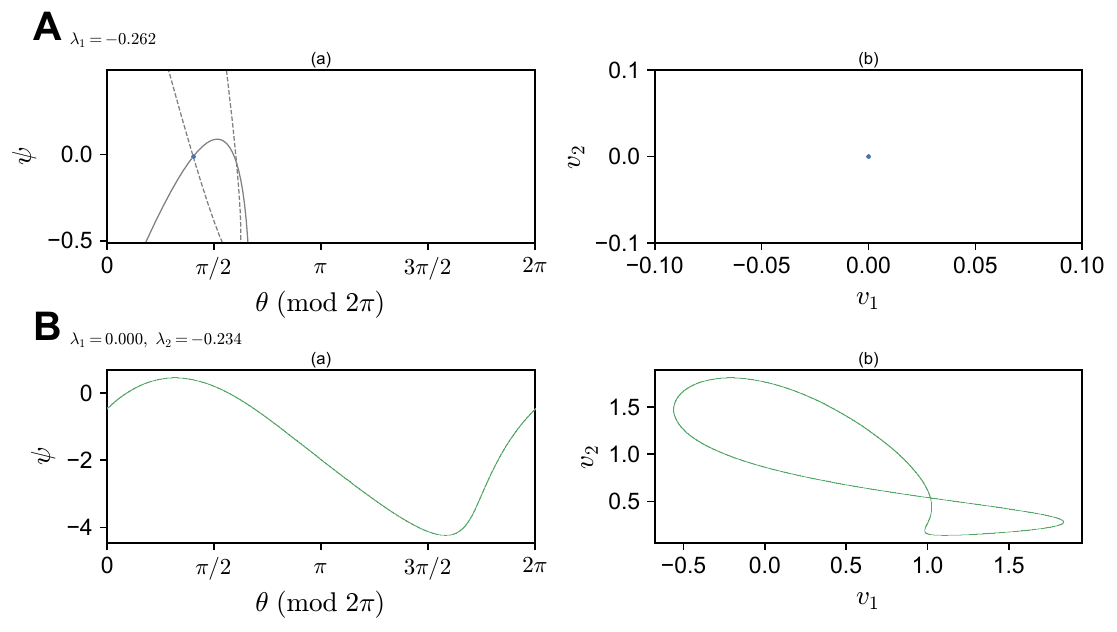}
  \caption{
    Representative cases A and B in the parameter map of \cref{fig:asym_le_spec_trajs}: (a) the $(\theta,\psi)$ projection and (b) the $(v_1,v_2)$ projection over the last $25{,}000$ time units with zero initial state.
    $\lambda_1$ and $\lambda_2$ show the estimates of two largest finite-time Lyapunov exponents using Benettin-QR algorithm.
  }
  \label{fig:asym_le_traj_ab}
\end{figure}

\begin{figure}[htbp]
  \centering
  \includegraphics[width=0.8\textwidth]{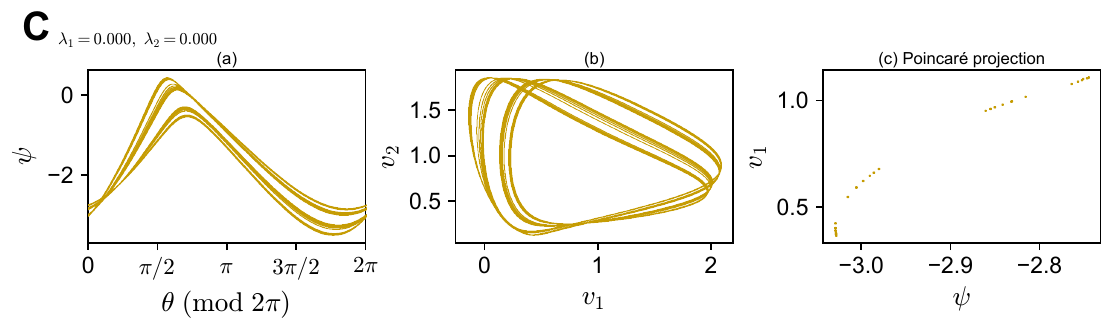}
  \caption{
    Representative case C in \cref{fig:asym_le_spec_trajs}. 
    Panels (a) and (b) give the $(\theta,\psi)$ and $(v_1,v_2)$ projections over the last $25{,}000$ time units with zero initial state; (c) plots the $(\psi,v_1)$ projection of the Poincar\'e section $\theta=0\pmod{2\pi}$ with $\dot\theta>0$, attained by linearly interpolating the saved trajectory. 
    $\lambda_1$ and $\lambda_2$ show the estimates of two largest finite-time Lyapunov exponents using Benettin-QR algorithm.
  }
  \label{fig:asym_le_traj_c}
\end{figure}

\begin{figure}[htbp]
  \centering
  \includegraphics[width=0.8\textwidth]{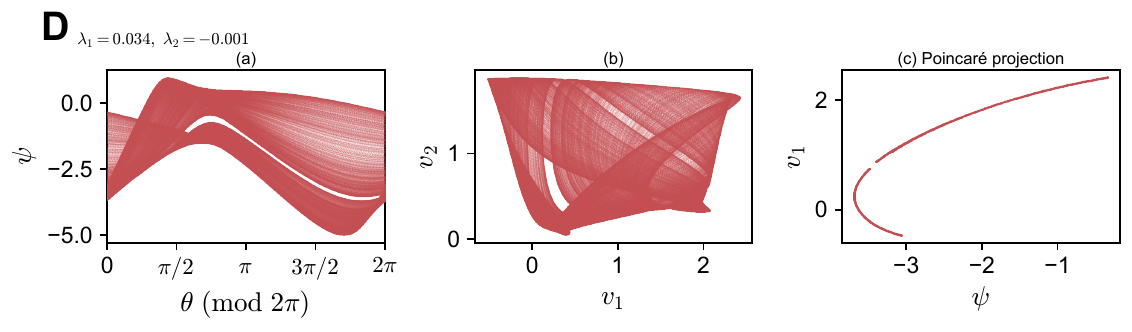}
  \caption{
    Representative case D in \cref{fig:asym_le_spec_trajs}. 
    Panel definitions and plotting settings are as specified in \cref{fig:asym_le_traj_c}.
  }
  \label{fig:asym_le_traj_d}
\end{figure}

\begin{figure}[htbp]
  \centering
  \includegraphics[width=0.8\textwidth]{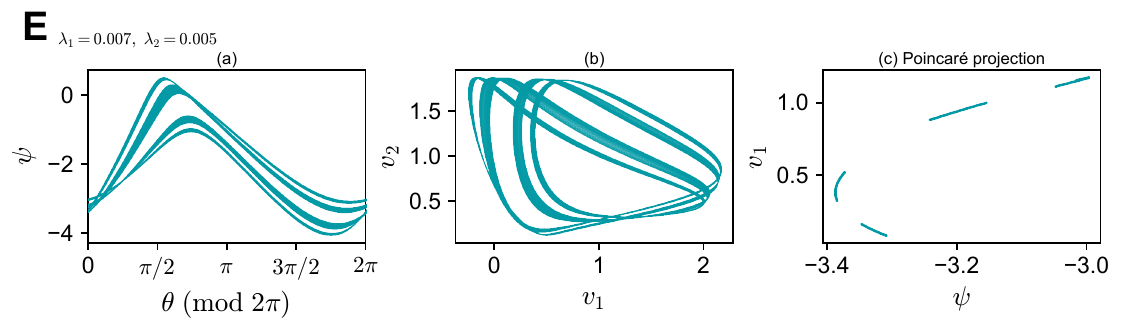}
  \caption{
    Representative case E in \cref{fig:asym_le_spec_trajs}. 
    The plotting setup stays the same as \cref{fig:asym_le_traj_c}, with $\theta$ displayed modulo $2\pi$ in panel (a).}
  \label{fig:asym_le_traj_e}
\end{figure}

For the asymmetric map in \cref{fig:asym_le_spec_trajs}, we set $\beta_1=0.487$, $\beta_2=0.929$, $\kappa_2=0.08$ and $i_2=0.95$ with initial state $(\theta, \psi, v_1, v_2) = (0, 0, 0, 0)$, and let $\kappa_1 \in [0.1, 0.6]$ and $i_1 \in [0.9, 1.5]$ vary. 
We set a transient time $T_{\mathrm{tr}}=300$ units and measurement time $T_{\mathrm{me}}=1500$ units. 
Blue, green, yellow, red and cyan regions correspond to spectral signatures consistent with stable fixed points, limit cycles, quasi-periodic tori, chaos and hyperchaos respectively, reported by the finite-time Lyapunov exponents. 
Unlike the symmetric parameter scan in \cref{fig:sym_le_spec_trajs}, all five spectral classes occur at different parameter values in this asymmetric parameter plane, with a mixed chaotic and hyperchaotic region that grows around $(i_1=1.04, \kappa_1=0.24)$.
Moreover, the chaotic region is separated by strips composed of limit cycles and quasi-periodic tori, forming a flame-like structure. 
We also plot the projected dynamics of the five selected parameter settings (A)--(E) and the projected Poincar\'e sections of (C)--(E) over the last $25,000$ of $50,000$ time units of simulation in figures \ref{fig:asym_le_traj_ab}-\ref{fig:asym_le_traj_e}. 
The upper setting (A) ($(i_1,\kappa_1)=(0.96012024,\,0.4496994)$, $\lambda_1 = -0.262$) in \cref{fig:asym_le_traj_ab} belongs to the stable-fixed-point region in \cref{fig:asym_le_spec_trajs}.
The system settles into a stable fixed point at an intersection of the $v_1$ and $v_2$ nullclines shown in (a). 
In the lower setting (B) ($(i_1,\kappa_1)=(1.4002004,\,0.42965932)$, $\lambda_1 \approx 0, \lambda_2 = -0.234$) selected in \cref{fig:asym_le_traj_ab}, the system exhibits periodically oscillating dynamics in panel (a) and a closed curve in the $(v_1, v_2)$ projection, indicating a rotating limit cycle. 
In setting (C) ($(i_1,\kappa_1)=(1.19458918,\,0.30240481)$, $\lambda_1, \lambda_2 \approx 0$) illustrated in \cref{fig:asym_le_traj_c}, the dynamics are consistent with quasi-periodic motion on a torus, whose projections on $(\theta, \psi)$ and $(v_1, v_2)$ exhibit modulated oscillations.
The projected Poincar\'e section of (C) has separated points, distinguishing it from the continuous curve segments in (D) and (E) sketched in \cref{fig:asym_le_traj_d} and \cref{fig:asym_le_traj_e}, respectively.
In setting (D) ($(i_1,\kappa_1)=(1.05511022,\,0.17715431)$, $\lambda_1 = 0.034, \lambda_2 = -0.001$) and E ($(i_1,\kappa_1)=(1.23426854,\,0.26833667)$, $\lambda_1 = 0.007, \lambda_2 = 0.005$), both projected trajectories in (a) and (b) evolve irregularly.
Additionally, the Poincar\'e sections of (D) and (E) have continuous curve segments, with (D) having two segments and (E) having four.

\begin{figure}[htbp]
  \centering
  \includegraphics[width=0.8\textwidth]{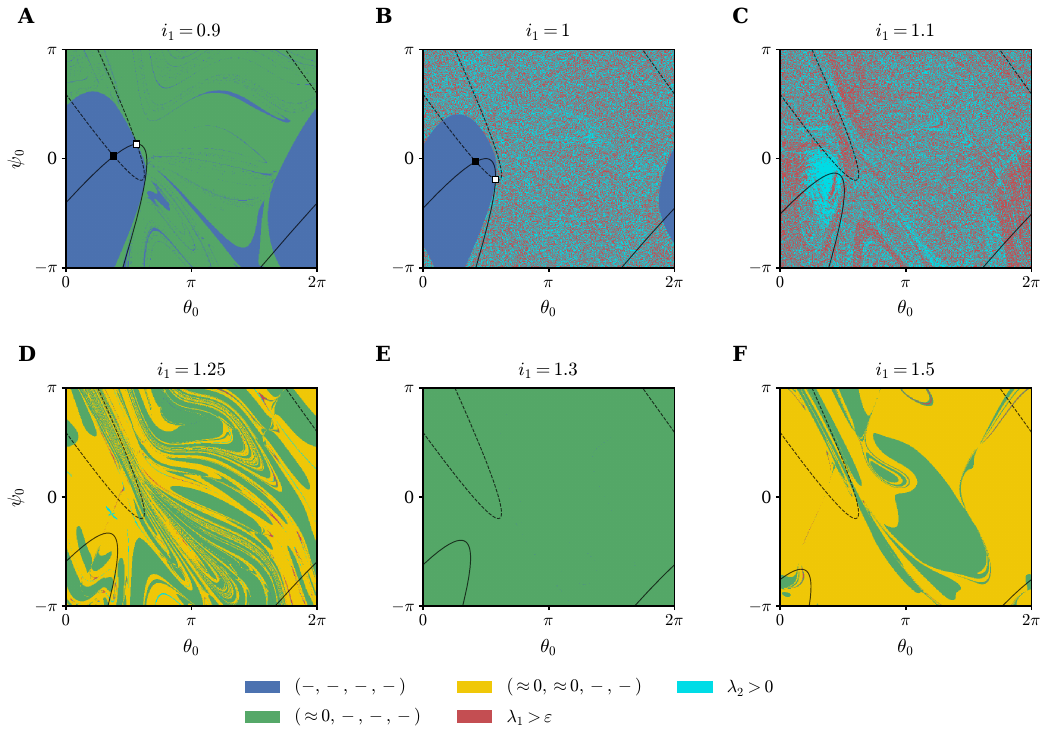}
  \caption{
    Finite-time deterministic Lyapunov classifications while varying asymmetric current $i_1$ at fixed $\kappa_1=0.24$ with $(\beta_1,\beta_2,i_2,\kappa_2)$ fixed at $(0.487,0.929,0.95,0.08)$. 
    The initial-condition setup, and nullcline markers are specified in \cref{fig:det_le_ic_0} and color encodings are stated in \cref{fig:asym_le_spec_trajs}.
  }
  \label{fig:det_le_ic_2}
\end{figure}

\begin{figure}[htbp]
  \centering
  \includegraphics[width=0.7\textwidth]{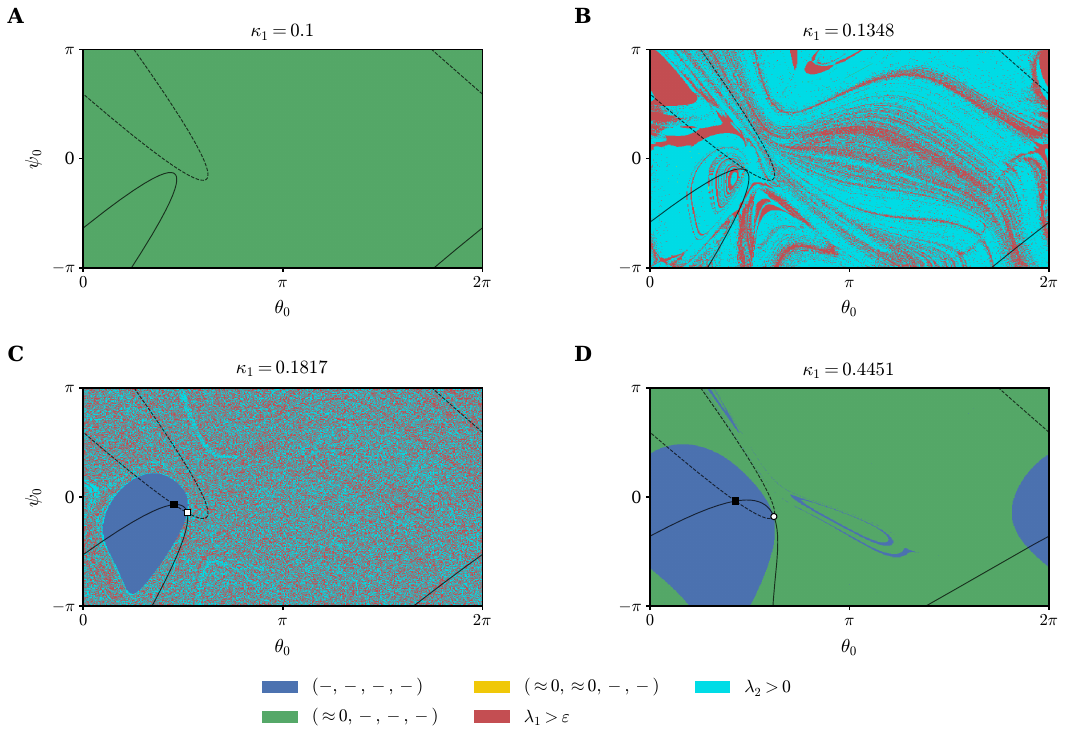}
  \caption{
    Finite-time deterministic Lyapunov classifications while varying asymmetric coupling $\kappa_1$ at fixed $i_1=1.04$ with $(\beta_1,\beta_2,i_2,\kappa_2)$ fixed at $(0.487,0.929,0.95,0.08)$.  
    The initial-condition sampling, colors and markers are those used in \cref{fig:det_le_ic_2}.
  }
  \label{fig:det_le_ic_3}
\end{figure}

The asymmetric initial-condition maps \cref{fig:det_le_ic_2,fig:det_le_ic_3} include six drive values, $i_1=0.9, 1.0, 1.1, 1.25, 1.3, 1.5$, and four coupling values, $\kappa_1=0.1, 0.1348, 0.1817, 0.4451$. 
Each panel in figures \ref{fig:det_le_ic_2}-\ref{fig:det_le_ic_3} uses a $500\times500$ grid in $(\theta_0,\psi_0)$, with $v_1(0)=v_2(0)=0$, $\theta_0\in[0,2\pi]$, and the displayed $\psi_0$ range. 
In all asymmetric panels in \cref{fig:det_le_ic_2} and \cref{fig:det_le_ic_3} except that at $i_1=1.25$, points initially classified as having two near-zero exponents or $\lambda_1>\varepsilon$ are recomputed with $T_{\mathrm{tr}}=300$ and $T_\mathrm{me}=1500$ for final classification. 

In \cref{fig:det_le_ic_2}, blue and green regions coexist at $i_1=0.9$, while red and cyan occupy most of the displayed plane at $i_1=1.0$ and $1.1$.
At $i_1=1.25$, green and yellow bands fill the plane.
The plane is almost entirely green at $i_1=1.3$, whereas yellow predominates again at $i_1=1.5$.
In \cref{fig:det_le_ic_3}, the map is green at $\kappa_1=0.1$ while it becomes predominantly red and cyan at $\kappa_1=0.1348$ and $0.1817$.
At $\kappa_1=0.4451$, this map contains blue and green regions when the $v_1$ and $v_2$ nullclines intersect again at larger $\kappa_1$.

\section{Bifurcation Diagram}\label{sec:bifurcation_diagram}

\begin{figure}[htbp]
  \centering
  \includegraphics[width=0.9\textwidth]{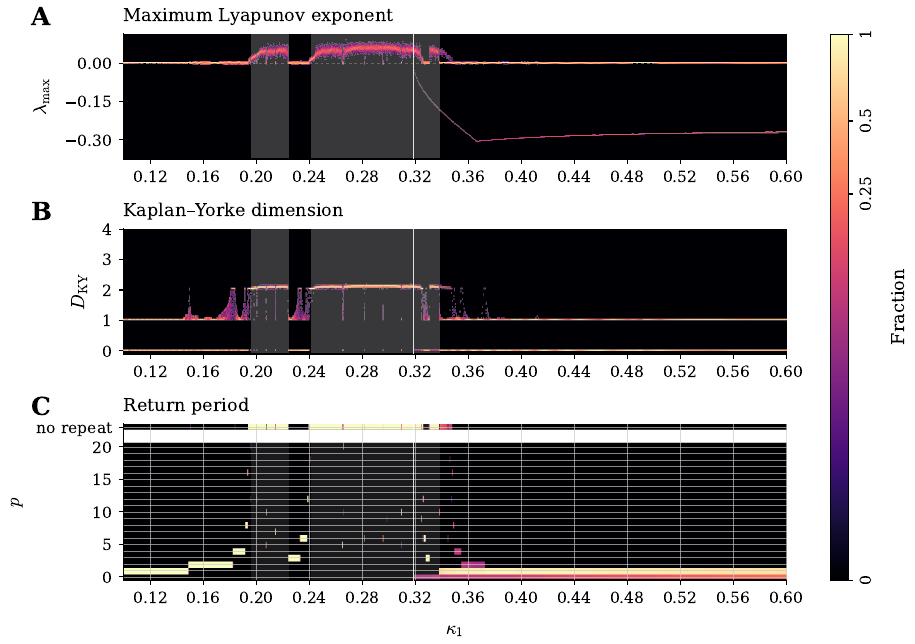}
  \caption{Finite-window deterministic scan at $i_1=1.1056$, varying $\kappa_1\in[0.1,0.6]$, with $(\beta_1,\beta_2,i_2,\kappa_2)=(0.487,0.929,0.95,0.08)$. Each of $1{,}000$ parameter values uses $120$ equally weighted initial conditions on a $10\times12$ grid in $(\theta_0,\psi_0)\in[0,2\pi)\times[-\pi,\pi]$ and $v_1(0)=v_2(0)=0$. 
  The deterministic Benettin--QR method uses RK4 step $0.02$ and QR interval $1$. 
  The observation window is $t\in[2000,2400]$. Panels A, B and C show the distributions of the maximum Lyapunov estimate, raw Kaplan--Yorke estimate and of period of Poincar\'e return to the section $\phi_2=2\pi n, n \in Z$ with increasing $\phi_2$. 
  The row $p=0$ denotes trajectories with no upward crossing in the observation window; ``no repeat'' denotes trajectories with no detected period of crossings.
  Colors give fractions of the $120$ initial conditions per bin, on a nonlinear scale. The state distributions are shown in \cref{fig:sym_buf_i11056_states}.
  }
  \label{fig:sym_buf_i11056}
\end{figure}

\begin{figure}[htbp]
  \centering
  \includegraphics[width=0.8\textwidth]{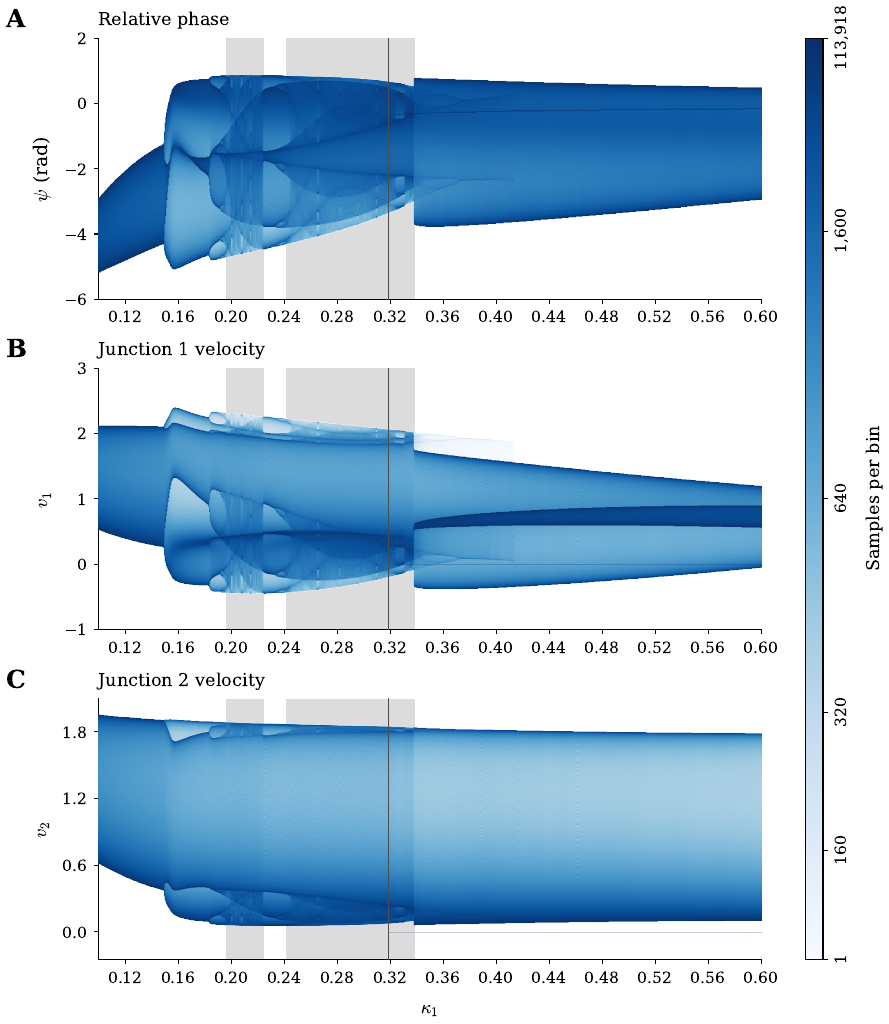}
  \caption{Finite-window deterministic scan at $i_1=1.1056$, varying $\kappa_1\in[0.1,0.6]$, with $(\beta_1,\beta_2,i_2,\kappa_2)=(0.487,0.929,0.95,0.08)$. Each of $1{,}000$ parameter values uses $120$ equally weighted initial conditions on a $10\times12$ grid in $(\theta_0,\psi_0)\in[0,2\pi)\times[-\pi,\pi]$ and $v_1(0)=v_2(0)=0$. The deterministic Benettin--QR method uses RK4 step $0.02$ and QR interval $1$. The observation window is $t\in[2000,2400]$. Panels A--C show time-sample counts for $\psi$, $v_1$, and $v_2$. Darker blue indicates larger counts, with $480{,}120$ samples per parameter and a nonlinear color scale. The corresponding spectral distributions are shown in \cref{fig:sym_buf_i11056}.}
  \label{fig:sym_buf_i11056_states}
\end{figure}

\begin{figure}[htbp]
  \centering
  \includegraphics[width=0.9\textwidth]{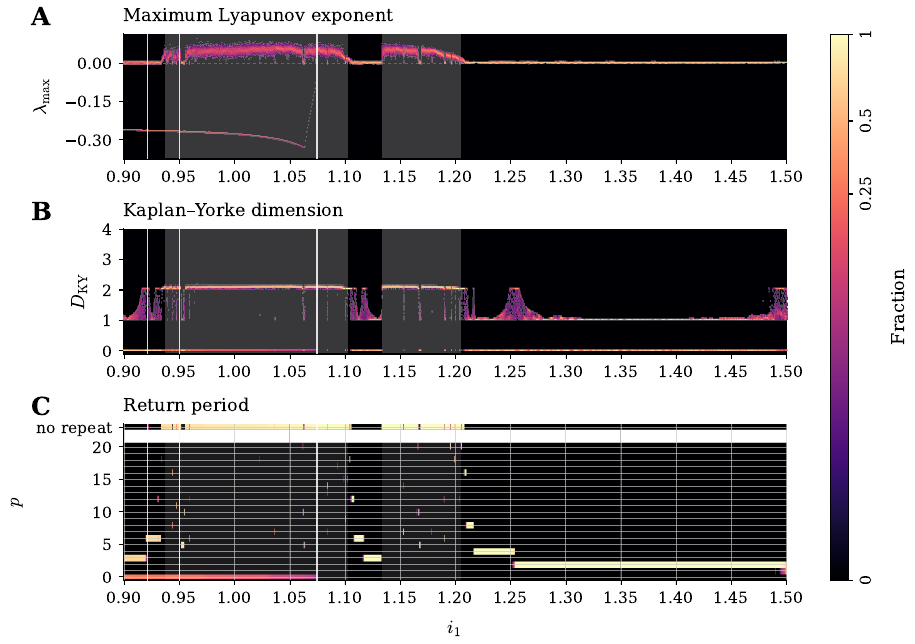}
  \caption{
    Finite-window deterministic scan of MFLE, KY dimension and Poincar\'e return period at $\kappa_1=0.24$, varying $i_1\in[0.9,1.5]$, with $(\beta_1,\beta_2,i_2,\kappa_2)=(0.487,0.929,0.95,0.08)$.  
    The state distributions are shown in \cref{fig:asym_buf_i1_k24_states}.
    All plotting conventions are the same as in \cref{fig:sym_buf_i11056}.
  }
  \label{fig:asym_buf_i1_k24}
\end{figure}

\begin{figure}[htbp]
  \centering
  \includegraphics[width=0.8\textwidth]{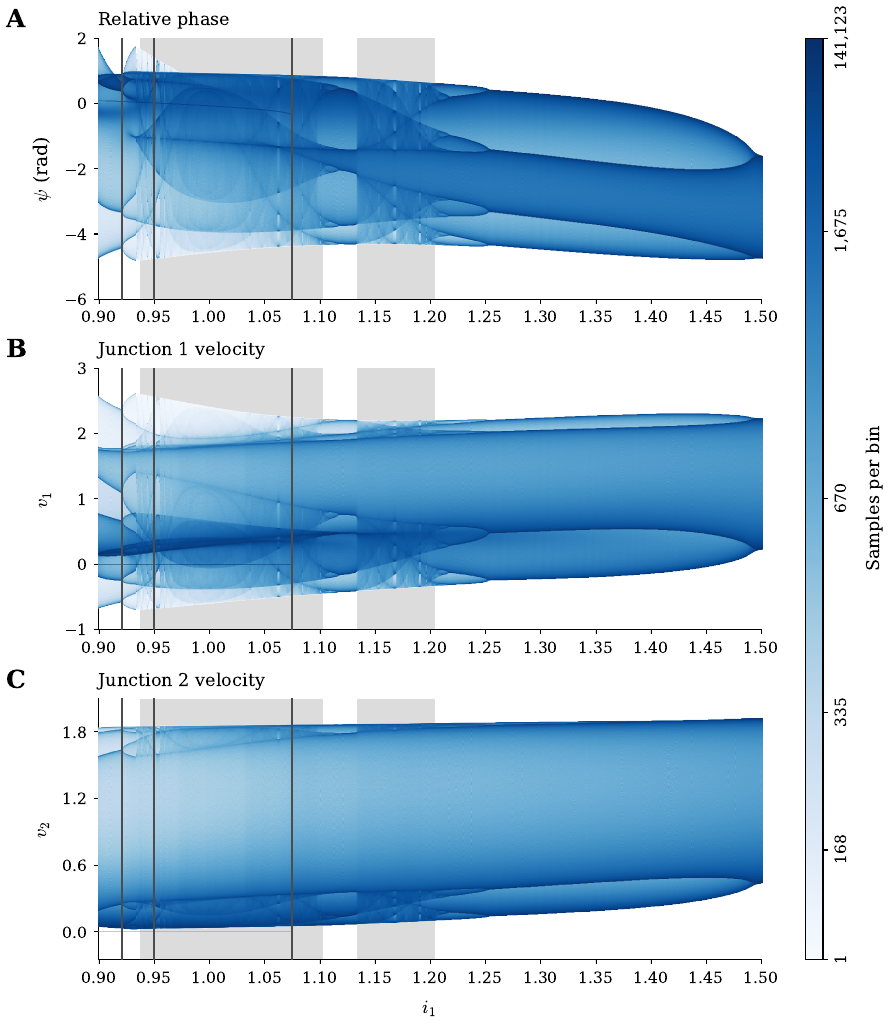}
  \caption{
    Finite-window deterministic scan of the state distributions $(\psi,v_1,v_2)$ at $\kappa_1=0.24$, varying $i_1\in[0.9,1.5]$, with $(\beta_1,\beta_2,i_2,\kappa_2)=(0.487,0.929,0.95,0.08)$. 
    The corresponding spectral distributions are illustrated in \cref{fig:asym_buf_i1_k24}.
    Sampling and display settings follow \cref{fig:sym_buf_i11056_states}.
  }
  \label{fig:asym_buf_i1_k24_states}
\end{figure}

\begin{figure}[htbp]
  \centering
  \includegraphics[width=0.9\textwidth]{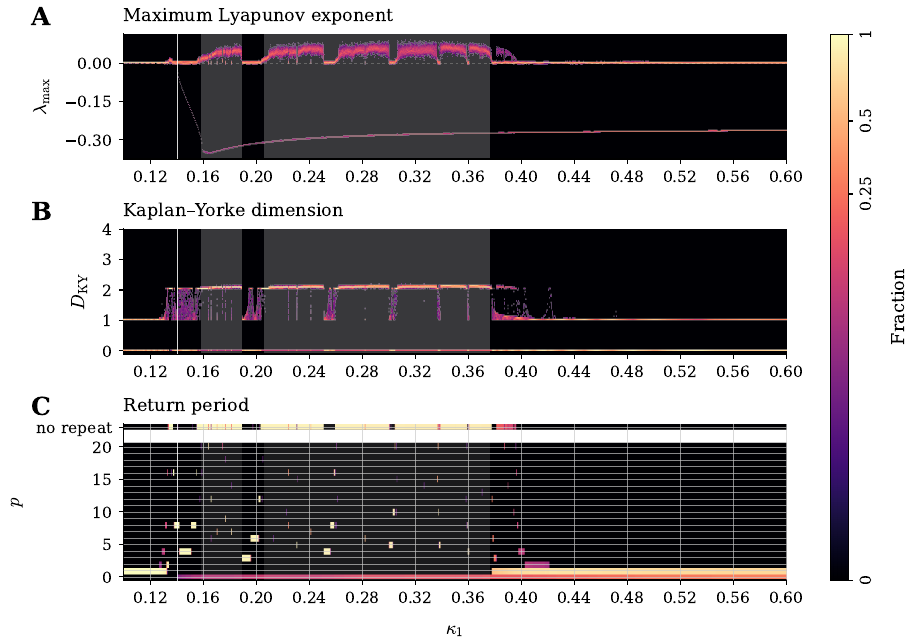}
  \caption{
    Finite-window deterministic scan of MFLE, KY dimension and Poincar\'e return period at $i_1=1.04$, varying $\kappa_1\in[0.1,0.6]$, with $(\beta_1,\beta_2,i_2,\kappa_2)=(0.487,0.929,0.95,0.08)$. 
    The associated state distributions are presented in \cref{fig:asym_buf_k1_i104_states}. 
    The sampling and display conventions are the same as in \cref{fig:sym_buf_i11056}.
  }
  \label{fig:asym_buf_k1_i104}
\end{figure}

\begin{figure}[htbp]
  \centering
  \includegraphics[width=0.8\textwidth]{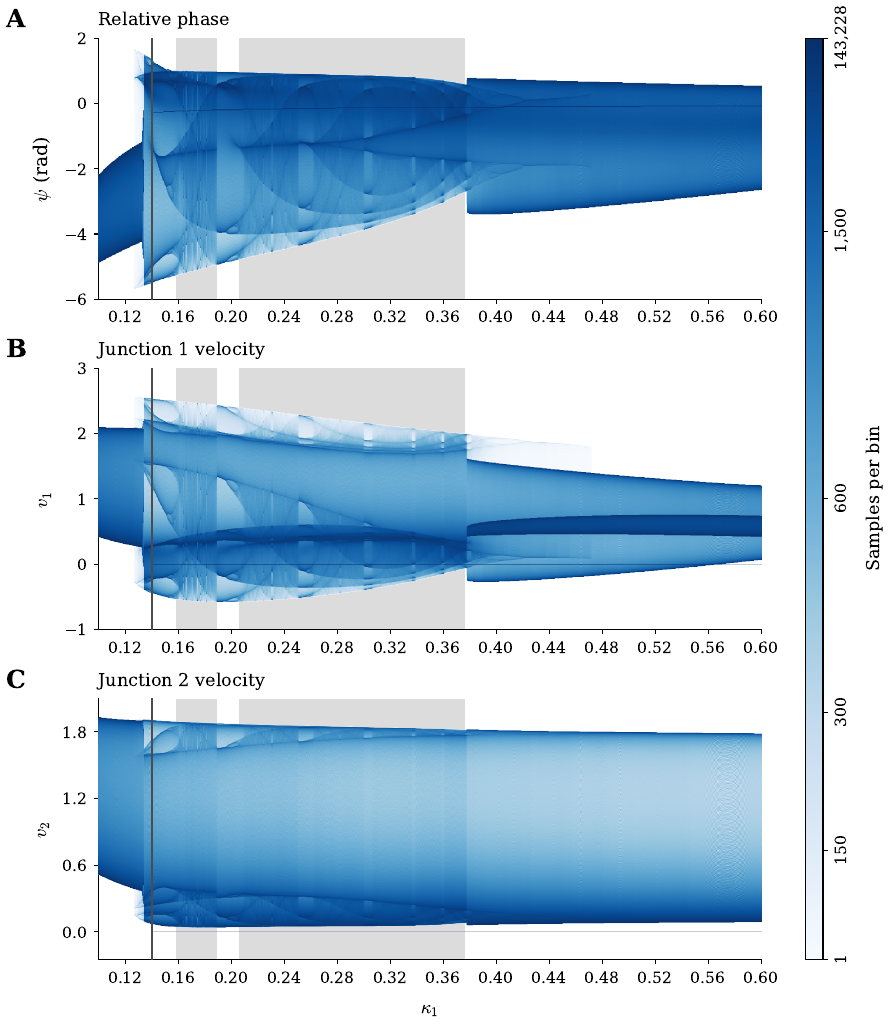}
  \caption{
    Finite-window deterministic scan of the state distributions $(\psi,v_1,v_2)$ at $i_1=1.04$, varying $\kappa_1\in[0.1,0.6]$, with $(\beta_1,\beta_2,i_2,\kappa_2)=(0.487,0.929,0.95,0.08)$. 
    The corresponding spectral distributions are plotted in \cref{fig:asym_buf_k1_i104}.
    Sampling, panel layout and color scales are stated in \cref{fig:sym_buf_i11056_states}.
    }
  \label{fig:asym_buf_k1_i104_states}
\end{figure}

To resolve the dependence on both parameters and initial conditions, figures \ref{fig:sym_buf_i11056}-\ref{fig:asym_buf_k1_i104_states} combine the sampled distributions of $(\psi,v_1,v_2)$ with distributions of the maximum finite-time Lyapunov exponents (MFLE) and the Kaplan--Yorke (KY) estimate across $120$ equally weighted initial conditions.
Each scan contains $1{,}000$ parameter values, and all observations use the window $t\in[2000,2400]$ after discarding the preceding transient dynamics.
The spectra are computed by the Benettin--QR method using double-precision RK4 with step $0.02$ and QR orthogonalization every unit time.
Figures \ref{fig:sym_le_spec_trajs}--\ref{fig:det_le_ic_3} in the previous section describe variation across parameter settings as well as initial conditions, whereas the spectral panels (\ref{fig:sym_buf_i11056}, \ref{fig:asym_buf_i1_k24} and \ref{fig:asym_buf_k1_i104}) summarize distributions across the $120$ initial conditions, and state panels (\ref{fig:sym_buf_i11056_states}, \ref{fig:asym_buf_i1_k24_states} and \ref{fig:asym_buf_k1_i104_states}) here pool time samples from their trajectories.
The state distributions therefore describe the range explored by periodic, chaotic, or coexisting trajectories.
Writing the ordered finite-window spectrum as $\lambda_1\ge\cdots\ge\lambda_4$, the plotted KY dimension estimates are
\begin{equation}\label{eq:ky_scan_estimate}
  \widehat D_{\mathrm{KY}}(x_0) = 
  \begin{cases}
    0, & \lambda_1 < 0,\\
    4, & \sum_{\ell=1}^{4}\lambda_\ell(x_0) \ge 0, \\
    \displaystyle j+\frac{\sum_{\ell=1}^{j}\lambda_\ell(x_0)}{|\lambda_{j+1}(x_0)|}, & \mathrm{otherwise},
  \end{cases}\qquad,
\end{equation}
where $1\le j\le3$ is the largest index with a nonnegative partial sum. 

To interpret representative parameter values, we also examine returns of the trajectories to the section $\phi_2=2\pi n$ with increasing $\phi_2$, recording $(\psi,v_1,v_2)$ at each crossing.
A return period $p$ means that these states approximately repeat after $p$ crossings, with each crossing obtained by linear interpolation.
The smallest $p\le20$ with at least $2p+1$ crossings satisfying a maximum return-state discrepancy below $5\times10^{-3}$ is used as a criterion of repetition.
We classify trajectories as approximately stationary when both velocity magnitudes remain below $10^{-5}$ throughout the saved window.

For the coupling scan at $i_1=1.1056$ in \cref{fig:sym_buf_i11056}, the low-coupling regime ($\kappa < 0.2$) is dominated by regular running with a maximum exponent close to zero.
Representative trajectories at $\kappa_1\approx0.14004$, $0.16006$, $0.19009$ and $0.19309$ have return periods $1$, $2$, $4$ and $8$ in panel (C), respectively, supporting a period-doubling route towards the first broad positive-exponent region in panel (A), located approximately at $0.195\lesssim\kappa_1\lesssim0.224$.
In \cref{fig:sym_buf_i11056_states}, the edges of occupancy get blurred in the two positive-exponent regions, which is most easily seen in panels (A) and (B).
Within this positive-exponent region containing $\kappa_1 = 0.2$ in panel (A) of \cref{fig:sym_buf_i11056}, the maximum exponents are mostly positive, and KY dimension estimates in panel (B) are slightly above $2$, providing evidence of chaotic running over the observation window.
Narrow downward excursions interrupt this pattern.
Representative trajectories near $\kappa_1=0.20811$ and $0.21461$ exhibit return periods $10$ and $7$ in panel (C) of \cref{fig:sym_buf_i11056}, respectively, supporting embedded periodic windows in the chaotic region illustrated in all panels in \cref{fig:sym_buf_i11056_states}.

A wider regular window appears before the second chaotic region shown in MFLE distribution of \cref{fig:sym_buf_i11056} at larger $\kappa_1$.
Near $\kappa_1 \approx 0.22513$ and $\kappa_1 \approx 0.23013$, representative trajectories yield yet another period doubling scheme in Poincar\'e return period spectrum of \cref{fig:sym_buf_i11056}. 
The return period of the system increases from $3$ to $6$ near $\kappa_1 \approx 0.23514$ and eventually reaches a period of $12$ at $\kappa_1 \approx 0.24$.
The second broad chaotic region reappears, starting from $\kappa_1 \approx 0.241$ to approximately $\kappa_1 \approx 0.325$, including narrow interruptions with period-$5$, period-$6$ and period-$10$.
A set of periodic windows with periods that are multiples of $3$ occurs near the small window at $\kappa_1 \approx 0.33023$, where the MFLE is near-zero and KY dimension is lower.
Positive exponents resume from $\kappa_1 \approx 0.331$ to $\kappa_1 \approx 0.338$.

Note that a distinct transition occurs in the MFLE panel within the second chaotic region triggered by a fixed-point bifurcation.
Solving the equilibrium equations \eqref{eq:coupled_jj-system-transformed}, we locate a saddle-node $(\theta_{*},\psi_{*}, v_{1,*}, v_{2,*})\approx(1.59,-0.41, 0, 0)$ at $\kappa_{1,\mathrm{SN}}\approx0.318637$ with eigenvalues approximately $0$, $-0.6994$, and $-0.3583\pm0.42877\mathrm{i}$ of the Jacobian in \eqref{Jacobian}. 
The zero eigenvalue of this Jacobian is simple. 
The parameter and quadratic coefficients of the local reduction are nonzero. 
The two nearby equilibrium branches on the larger-coupling side have zero and one unstable direction, respectively, creating a stable equilibrium and a saddle corresponding to a negative branch in MFLE.
At $\kappa_1\approx0.32022$, for instance, the stationary and running regimes coexist in the right halves of all projection panels of \cref{fig:sym_buf_i11056_states}: four of the $120$ saved trajectories are approximately stationary, whereas the other $116$ keep running.                             

The negative-exponent branch has a turn near $\kappa_1=0.36$--$0.37$.
When $\kappa_1=0.36$, the stable focus $(\theta_{*},\psi_{*}, v_{1,*}, v_{2,*})\approx(1.493,-0.301, 0, 0)$ has two real Jacobian eigenvalues $\mu_1=-0.277978$, $\mu_2=-0.522695$ and a complex-conjugate pair $\mu_{H,\pm}=-0.307663\pm0.560312\mathrm{i}$.
The real eigenvalue $\mu_1$ has the largest real part thus it determines the slowest decay of perturbations around the stationary state.
As $\kappa_1$ increases, $\mu_1$ decreases while the real part of $\mu_{H,\pm}$ increases, until $\mu_1 = \Re\,\mu_{H, \pm} \approx -0.304519$, at $\kappa_1\approx0.366051$.
At this point, $\mu_{H,\pm}\approx-0.304519\pm0.571515\mathrm{i}$ and $\mu_2\approx-0.502442$.
When $\kappa_1=0.37$, the eigenvalues become $\mu_1=-0.323213$, $\mu_2=-0.487474$ and $\mu_{H,\pm}=-0.302657\pm0.578527\mathrm{i}$, so the complex-conjugate pair now determines the slowest decay.
The turn of the negative branch in panel (A) of \cref{fig:sym_buf_i11056} therefore corresponds to a change in the dominant decay when a non-oscillatory mode switches to an oscillatory mode, with all eigenvalue real parts remaining negative.
The two real eigenvalues remain distinct through this turn and meet later at $\kappa_1\approx0.379939$, where $\mu_1=\mu_2\approx-0.409486$.
They then form a second complex-conjugate pair, while the original pair retains the largest real part.
Thus, the change in the dominant mode precedes the merger of the two real eigenvalues, and the stationary state remains stable through both changes.

Beyond approximately $\kappa_1=0.338$, the fraction of sampled initial conditions with positive exponents decreases, while regular running becomes more prominent alongside the stationary branch visible in panel (A) of \cref{fig:sym_buf_i11056}.
The counts of trajectories whose $\lambda_1 > 0.01$ are $111$, $20$, $12$, and $0$ at $\kappa_1\approx0.33524$, $0.34024$, $0.34525$, and $0.34975$, respectively.
The saved trajectories include low-period running states in this transition region, and at $\kappa_1=0.6$ there are $28$ approximately stationary trajectories and $92$ running trajectories.

The current scan at {$\kappa_1=0.24$} showcased in \cref{fig:asym_buf_i1_k24} starts with a stationary branch visible in the MFLE distribution panel.
This branch shrinks in the number of trajectory occupancies and disappears as the normalized current parameter increases until the fixed point bifurcation at $i_1 \approx 1.075$.
When $i_1 \le 0.92$, a stable stationary state already coexists with periodic running illustrated in \cref{fig:asym_buf_i1_k24_states}: $35$ trajectories are approximately stationary at $i_1 \approx 0.9$, while the remaining $85$ running trajectories have dominant return period $3$ as the pronounced bright stripe in panel (C) of \cref{fig:asym_buf_i1_k24}.
Period-$6$ returns are found between $i_1\approx0.92523$ and $i_1\approx0.93003$, followed by a thin window of return period $9$ and a broad positive-exponent region that begins at $i_1\approx0.935$.
The first chaotic region in panel (A) of \cref{fig:asym_buf_i1_k24} ends at $i_1 \approx 1.104$, and is interrupted by narrow regular windows containing the largest window of period $5$ in panel (C).
Equilibrium and Jacobian calculations locate an SN bifurcation at $i_{1,\mathrm{SN}}\approx1.075$, where a stable equilibrium and a saddle collide in a nondegenerate saddle-node bifurcation, corresponding to the rightmost vertical white line in both \cref{fig:asym_buf_i1_k24} and \cref{fig:asym_buf_i1_k24_states}.
The other two vertical lines in all panels of both \cref{fig:asym_buf_i1_k24,fig:asym_buf_i1_k24_states} correspond to two SS bifurcations, with one creating two saddles and one eliminating two saddles.
When $i_1 = 1.06$, the stable focus $(\theta_{*},\psi_{*}, v_{1,*}, v_{2,*})\approx(1.45, -0.25, 0, 0)$ has two sets of conjugate eigenvalues $\mu_{H, \pm} = -0.322\pm 0.499\mathrm{i}$ and $\mu_{L, \pm} = -0.386\pm 0.111\mathrm{i}$.
As $i_1$ increases, the pair $\mu_{L, \pm}$ merge at $i_1 \approx 1.063$ where $\mu_{L, \pm} \approx -0.379$, and then separate into two distinct real eigenvalues.
At $i_1=1.0635$ these are $\mu_{1} \approx -0.349$ and $\mu_{2} \approx -0.409$.
Following the splitting of the complex-conjugate pair into two real eigenvalues, the larger real eigenvalue overtakes the real part of the remaining complex pair, changing the dominant decay mode and causing the bend in the negative branch of MFLE in \cref{fig:asym_buf_i1_k24}.

The first broad chaotic region in MFLE distribution of \cref{fig:asym_buf_i1_k24} is followed by a regular interval from $i_1\approx1.105$ to $i_1 \approx 1.133$.
Representative return periods decrease from period-$18$ near $i_1 \approx 1.1048$ to period-$12$ near $i_1 \approx 1.1072$ to period-$6$ near $i_1\approx1.1102$. 
Period-$3$ returns occur between $i_1\approx1.1198$ and $i_1 \approx 1.13$. 

A second broad positive-exponent region then appears around $i_1=1.134$--$1.206$, again with narrow interruptions; for instance, period-$5$ returns occur near $i_1=1.16787$.
At its upper edge, higher-period rotation replaces the broad positive-exponent distribution in panel (A) of \cref{fig:asym_buf_i1_k24}.
Representative trajectories have return periods $16$ near $i_1\approx1.20991$, period $4$ near $i_1\approx1.22973$, period $2$ near $i_1\approx1.25976$ and $i_1\approx1.3$.

As $i_1$ increases, the regular intervals following the chaotic regions contain return-period subsequences with successive halvings: $12\rightarrow6\rightarrow3$ between the two broad chaotic regions and $16\rightarrow8\rightarrow4\rightarrow2$
beyond the second.
These subsequences resemble period-doubling cascades traversed in reverse.

These patterns also show up in \cref{fig:asym_buf_i1_k24_states}, where pairwise merging of the projected branches reduces the number of edges in the trajectory occupancy plots by half at each step.

At larger currents, a strong zero-KY band in panel (B) of \cref{fig:asym_buf_i1_k24} accompanies the sampled running motion, with a surge in KY dimension estimates caused by second-largest finite-time Lyapunov exponents located around $i_1 = 1.25$, as indicated in panel (B) of \cref{fig:asym_buf_i1_k24}.
The lower KY band in this part of panel B therefore contains running trajectories, which corresponds to rotating dynamics visible in all panels of \cref{fig:asym_buf_i1_k24_states}.
The phase and velocity distributions remain broad after the maximum exponent returns close to zero, describing the range explored by these trajectories, corresponding to clear projected edges in \cref{fig:asym_buf_i1_k24_states}.

For the coupling scan at the lower normalized current {$i_1=1.04$} plotted in \cref{fig:asym_buf_k1_i104}, representative trajectories have return period of $1$ for $\kappa_1\lesssim0.13$.
Figures \ref{fig:asym_buf_k1_i104}-\ref{fig:asym_buf_k1_i104_states} establish that chaotic and period doubling regions alternate four times.
Before the first chaotic region, there is a period doubling sequence following $1 \rightarrow 2 \rightarrow 4 \rightarrow 8 \rightarrow 16$ spotted in panel (C) of \cref{fig:asym_buf_k1_i104}, during which there is another period doubling occupying smaller range of $\kappa_1$, leading to a short positive MFLE region.

The stable equilibrium branch is created in a nondegenerate saddle-node at $\kappa_{1,\mathrm{SN}}\approx0.140163$, substantially below the value $0.318637$ when $i_1\approx1.1056$, depicted in the MFLE panel.
This branch exists prior to the first broad chaotic region at approximately $\kappa_1\approx0.157$--$0.189$, and remains available throughout the subsequent running-state transitions.
Note that there is a much less sharp turn of the negative MFLE branch around $\kappa_1 \approx 0.16$ compared to the sharper turns in \cref{fig:sym_buf_i11056,fig:asym_buf_i1_k24}.
When $\kappa_1 = 0.16$, the fixed point $(\theta_{*},\psi_{*}, v_{1,*}, v_{2,*}) \approx \allowbreak(1.45, -0.25, 0, 0)$ corresponds to the Jacobian
\[
J_{0.16}\approx
\begin{bmatrix}
0 & 0 & 0.5 & 0.5\\
0 & 0 & -1 & 1\\
0.004520 & 0.157740 & -0.487 & 0\\
-0.243084 & -0.201542 & 0 & -0.929
\end{bmatrix}
\]
with eigenvalues $\mu_{H, \pm} = -0.347\pm0.029\mathrm{i}$ and $\mu_{L, \pm} = -0.361 \pm 0.422\mathrm{i}$.
When $\kappa_1 \approx 0.17$ after the transition point, the stable focus moves to $(\theta_{*},\psi_{*}, v_{1,*}, v_{2,*}) \approx (1.439, -0.236, 0, 0)$, corresponding to the following Jacobian
\[
J_{0.17}\approx
\begin{bmatrix}
0 & 0 & 0.5 & 0.5\\
0 & 0 & -1 & 1\\
-0.014272 & 0.177136 & -0.487 & 0\\
-0.247564 & -0.203782 & 0 & -0.929
\end{bmatrix}.
\] 
Dominance switches from the $\mu_{H, \pm}$ pair to the $\mu_{L, \pm}$ pair as $\kappa_1$ increases from $0.16$ to $0.17$. 

Near $\kappa_1=0.378$, the broad positive-exponent distribution gives way to low-period running, and the state-variable envelopes reorganize. 
A nested period cascading pattern takes place between $\kappa_1 = 0.378$ and approximately $\kappa = 0.42$, where the chaotic windows occupy by small fraction of the window.
The $v_1$ distributions in \cref{fig:asym_buf_k1_i104_states} show the same cascade in their uppermost branch. 
The occupancy reveals a sequence of collisions: the edges blur in the chaotic region, and then sharpen again while the branches collide pairwise, until a single prominent branch remains.    

Among KY dimension distributions in figures \ref{fig:sym_buf_i11056}, \ref{fig:asym_buf_i1_k24} and \ref{fig:asym_buf_k1_i104}, we find the estimated KY dimension sensitive to weakening transverse contraction. 
The distributions of KY dimension estimates spread upwards before a conspicuous increase of the maximum exponent.
When $\lambda_1 \ge 0$ and $\lambda_1 + \lambda_2 < 0$, the definition of the KY dimension estimate yields $\widehat D_{\mathrm{KY}}=1+\lambda_1/|\lambda_2|$, which increases as $\lambda_2$ approaches zero from below.
When $\lambda_1 + \lambda_2 \ge 0$ but $\lambda_1 + \lambda_2 + \lambda_3 < 0$, we have $\widehat D_{\mathrm{KY}} = 2+(\lambda_1+\lambda_2)/|\lambda_3| \ge 2$.
For one sampled initial condition as in panel (B) of \cref{fig:sym_buf_i11056}, the leading pair changes from approximately $(\lambda_1, \lambda_2) = (0.002576,-0.12334)$ at $\kappa_1\approx0.16006$ to $(\lambda_1, \lambda_2) = (0.00059,-0.000393)$ at $\kappa_1\approx0.18208$, while the KY dimension estimates rise from $1.0209$ to $2.0004$.
Larger KY dimension estimates therefore reflect the weakening of transverse contraction, as the small rise near $\kappa_1\approx0.149$ and the more pronounced rises near $\kappa_1\approx0.182$, $\kappa_1\approx0.233$, and $\kappa_1\approx0.238$ coincide with the return period doublings.

\section{Conclusion and Future Challenges}\label{sec:conclusion}

In this work, we examine the deterministic dynamics of an inductively coupled RCSJ pair at zero external magnetic flux. 
Fixed-point bifurcation figures \ref{fig:sym_i_nullclines}, \ref{fig:sym_kappa_nullclines}, \ref{fig:asym_i2_nullclines} and \ref{fig:asym_kappa1_nullclines} illustrate the four-fixed-point group structures and reveal periodicity of successive pairs of SS and SN bifurcations around $\theta$.
Figures \ref{fig:num_fix_sym}--\ref{fig:num_fix_asym} reveal non-monotonic changes in the number of equilibria along coupling slices, which arise from the unequal deformation of the nullclines.
 
Projected bifurcation diagrams (figures \ref{fig:sym_buf_i11056_states}, \ref{fig:asym_buf_i1_k24_states} and \ref{fig:asym_buf_k1_i104_states}) produced by multiple-initial-condition scans establish intervals of coexistence between stable stationary states and running motion.
Finite-time Lyapunov spectra and Poincaré returns further resolve period-doubling sequences and periodic windows in the running dynamics. 
These results connect equilibrium structure, running-state transitions, and initial-condition-dependent behavior within the parameter regimes studied.
Several directions for further investigation of coupled junctions are outlined below.

During normal operation, a SQUID carries an external magnetic flux, which may include a time-varying component. 
This flux is the quantity the device measures~\cite{lefevre1992thermal}.
Taking external flux into consideration allows the dependence of equilibrium and running states, switching statistics, and voltage–flux characteristics on the applied flux to be examined.
Moreover, inferring system parameters given characteristic curves under known external magnetic flux changes could be an effective way to calibrate manufactured DC-SQUIDs, which turns out to be an inverse problem.
Physics-informed or simulation-based learning offers a candidate approach to parameter inference, with identifiability and uncertainty as central evaluation criteria.
Wiener processes are able to model thermal fluctuations in the resistor of the RCSJ model.
For a single resistively shunted junction, Yuan~\cite{YUAN2025116917} demonstrated a clockwise hysteresis cycle and a planar limit-cycle bifurcation in the deterministic current--voltage characteristic, and analyzed the stochastic dynamics and the probability density under thermal fluctuations through the Fokker--Planck equation.
Jia et al. \cite{Jia2003SymmetricDCSQUID} considered symmetric overdamped DC-SQUIDs with thermal noise and derived the Fokker--Planck equations of the proposed stochastic system before examining analytical formulae for SQUID characteristics, "ripple" phenomena and impact of thermal noise on the currents.
Detailed analysis of stationary distributions of asymmetric junction settings, varying noise strengths and stochastic bifurcations shall be further studied.
Extending this analysis to larger junction arrays could reveal emergent synchronization beyond pairwise coupling. The growth in the number of parameters then motivates population-level descriptions of larger arrays.
Parameter distributions, such as $i_k\sim p(i_k)$, provide a compact representation of junction heterogeneity.
Graph-based approaches offer a way to represent the array topology, with graph neural networks as candidate tools for exploring collective dynamics.

\bigskip
\noindent\textbf{Data Availability}

The numerical algorithms and source code that support the findings of this study are available from the corresponding author upon reasonable request.

\bigskip
\noindent\textbf{Declaration of competing interest}

No author associated with this paper has disclosed any potential or pertinent conflicts which may be perceived to have impending conflict with this work.

\bigskip
\noindent\textbf{Acknowledgements}

The authors thank Pak Hang Chris Lau for insightful discussions on stochastic mechanisms. This work was supported by the Guangdong Basic and Applied Basic Research Foundation (Grant No. 2025A1515012560) and the Guangdong Introduction Program (Grant No. 2023QN10X753).

\end{document}